\documentclass[11pt]{amsart}
\usepackage[a4paper,
            bindingoffset=0in,
            left=1in,
            right=1in,
            top=1.2in,
            bottom=1.1in,
            footskip=.25in]{geometry}

\usepackage{mathpazo}
\usepackage{euler}

\usepackage{setspace}
\usepackage[dvipsnames]{xcolor}

\usepackage{tikz}
\usepackage{tikz-cd}
\usepackage{caption}

\usepackage{mathtools}

\usepackage{graphicx}
\usepackage{enumerate} 

\usepackage{latexsym}
\usepackage{amssymb,hyperref}
\hypersetup{colorlinks,linkcolor={MidnightBlue},citecolor={Maroon},urlcolor={PineGreen}}
\usepackage[cp850]{inputenc}
\usepackage{epsfig}
\usepackage{amsthm}
\usepackage{amscd}
\usepackage{amsmath}
\usepackage{amsfonts}
\usepackage{graphics}
\newtheorem{theorem}{Theorem}[section]
\newtheorem{lemma}[theorem]{Lemma}
\newtheorem{corollary}[theorem]{Corollary} 
\newtheorem{prop}[theorem]{Proposition}

\newtheorem*{theorem*}{Theorem} 
\newtheorem*{corollary*}{Corollary}
\theoremstyle{definition}
\newtheorem{example}[theorem]{Example}
\newtheorem{remark}[theorem]{Remark}
\newtheorem{definition}[theorem]{Definition}
\newtheorem*{remark*}{Remark}

\newtheorem*{definition*}{Definition}
\newtheorem*{example*}{Example}
\newtheorem{conj}[theorem]{Conjecture}
\newtheorem{construction}[theorem]{Construction}

\usepackage{quiver}

\newtheoremstyle{named}{}{}{\itshape}{}{\bfseries}{}{.0em}{\thmnote{#3}}
\theoremstyle{named}
\newtheorem*{named}{}

\newcommand{\BC}{\mathbb C} 
\newcommand{\BR}{\mathbb R} 
 \newcommand{\BQ}{\mathbb Q}
 \newcommand{\BZ}{\mathbb Z}
\newcommand{\BF}{\mathbb F} 
 
\newcommand{\BP}{\mathbb P} \newcommand{\BG}{\mathbb G}

\newcommand{\CC}{\mathcal C} 
\newcommand{\CE}{\mathcal E} \newcommand{\CF}{\mathcal F}
\newcommand{\CG}{\mathcal G} \newcommand{\CH}{\mathcal H}
 
 \newcommand{\CL}{\mathcal L}
\newcommand{\CM}{\mathcal M} 
\newcommand{\CO}{\mathcal O} \newcommand{\CP}{\mathcal P}
 \newcommand{\CR}{\mathcal R}
 
\newcommand{\CU}{\mathcal U} 
 
\newcommand{\CY}{\mathcal Y} \newcommand{\CZ}{\mathcal Z}

\newcommand\smvee{\raise0.9ex\hbox{$\scriptscriptstyle\vee$}}
\newcommand{\Eff}{\overline{\text{Eff}}}
\newcommand{\Nef}{\text{Nef}}

\newcommand{\CHom}{\mathcal{H}om} 
\newcommand{\CExt}{\mathcal{E}xt}

\DeclareMathOperator{\Hom}{Hom} 
\DeclareMathOperator{\Ext}{Ext}

\DeclareMathOperator{\Mor}{Mor}

\DeclareMathOperator{\Spec}{Spec}

\DeclareMathOperator{\rank}{rank}

\DeclareMathOperator{\Pic}{Pic}

\DeclareMathOperator{\codim}{codim}

\newcommand{\comment}[1]{}

\title[Irreducibility of $\widetilde{M}_{1,0}(X)$]{Moduli Space of Genus One Curves on Cubic Threefold}
\author{Enhao Feng}
\address{Boston College, Chestnut Hill, 02467, MA, USA}
\email{fenge@bc.edu}

\begin{document}

\setstretch{1.1}

\begin{abstract}
Let $X$ be a smooth cubic threefold. By invoking ideas from Geometric Manin's Conjecture, we give a complete description of the main components of the Kontsevich moduli space of genus one stable maps $\overline{M}_{1,0}(X)$. In particular, we show that for degree $e\geqslant 5$, there are exactly two irreducible main components, of which one generically parametrizes free curves birational onto their images, and the other corresponds to degree $e$ covers of lines. As a corollary, we classify components of the morphism space $\Mor(E,X)$ for a general smooth genus one curve $E$.
\end{abstract}

\maketitle

\section{Introduction}

Let $X$ be a smooth projective variety over $\BC$. A natural object to study is the Kontsevich moduli space $\overline{M}_{g,0}(X)$ of genus $g$ stable maps on $X$. Knowledge of the irreducibility and dimension of its components usually sheds light on the arithmetic and the enumerative geometry of $X$. While many works have focused on the case of rational curves, i.e. $g = 0$, fewer have explored the land of higher genus curves due to the more complicated structure of this moduli space. In this paper, we prove an irreducibility statement for the space of genus one curves when $X$ is a smooth cubic threefold, analogous to the results of \cite{Sta00} for rational curves. In particular, we focus on the main components of $\overline{M}_{1,0}(X)$, i.e. the components which generically parametrize stable maps with smooth domains. Let $H$ be the generator of the Picard group of $X$. We denote the union of the main components of $\overline{M}_{1,0}(X)$ by $\widetilde{M}_{1,0}(X)$, and we let $\widetilde{M}_{1,0}(X,e)$ be the sublocus of $\widetilde{M}_{1,0}(X)$ parametrizing $H$-degree $e$ stable maps.

\begin{theorem}
    \label{IrredM10} Let $X$ be a smooth cubic threefold with $\Pic(X) = \BZ H$. For $H$-degree $e\geqslant 3$ and $e \neq 4$, the locus $\widetilde{M}_{1,0}(X,e)$ consists of exactly two irreducible components $R_e$ and $N_e$, where
    \begin{itemize}
        \item $R_e$ has dimension $2e$ and generically parametrizes free genus one curves which map birationally onto their images, and
        \item $N_e$ has dimension $2e+2$ and generically parametrizes degree $e$ covers of lines.
    \end{itemize} 
    When $e = 4$, the space $\widetilde{M}_{1,0}(X,4)$ consists of $R_4$, $N_4$, and an extra component parametrizing double covers of conics. 
\end{theorem} 

Using this theorem, we can deduce the irreducibility of the components parametrizing free curves in the morphism space:

\begin{theorem} 
    \label{IrredMor} Let $X$ be a smooth cubic threefold with $\Pic(X) = \BZ H$. Let $E$ be a smooth connected genus one curve general in moduli. For $e\geqslant 6$, there is a unique irreducible component of dimension $2e$ in $\Mor(E,X)$ parametrizing free curves of $H$-degree $e$. 
\end{theorem}

One of the main new techniques we employ in the proof of Theorem \ref{IrredM10} is the following result, which we refer to as the Movable Bend-and-Break.

\begin{theorem}[Movable Bend-and-Break] \label{Thm:MBB} 
    Let $X$ be a smooth cubic threefold. For $e\geqslant 5$, any free component $M\subset\widetilde{M}_{1,0}(X,e)$ contains a stable map $f: Z = Z_1 + Z_2 \to X$ satisfying the following conditions:
    \begin{itemize}
        \item $Z_1$ is a smooth genus one curve and $f|_{Z_1}$ is free.
        \item $Z_2$ is a rational curve and $f|_{Z_2}$ is free.
        \item $f|_{Z_i}$'s are birational onto their images.
    \end{itemize}
\end{theorem}

We describe the strategy to deduce Theorem \ref{IrredM10} using this technique in Section \ref{resandstra}, and we briefly illustrate the contrast against the $g=0$ case in Section \ref{comparisongenus0}.

\subsection{Geometric Manin's Conjecture}
Our approach to the problem is motivated by Geometric Manin's Conjecture, proposed in \cite{LT19a}, and further developed in a series of works including \cite{LST22, BLRT22, LRT25}. The conjecture is based on an influential heuristic discovered by Batyrev \cite{Bat88} that leads to the formulation of the original Manin's Conjecture over number fields. Let $C$ be a smooth connected genus $g$ curve, $X$ be a smooth Fano variety over $\BF_q$, and $\Mor(C, X)$ be the space of morphisms from $C$ to $X$. Batyrev's idea can be summarized in the following statements.

\begin{enumerate}[\hspace{2mm}1.]
    \item (\textbf{Exceptional Set}) There is a proper closed set $Z\subset X$ such that the components of $\Mor(C,X)$ parametrizing curves lying outside $Z$ have the expected dimension.
    
    \item (\textbf{Uniqueness}) For each nef and integral algebraic equivalence curve class $\beta$, there is a unique irreducible component $M$ in $\Mor(C,X)$ parametrizing curves with class $\beta$.
    
    \item (\textbf{Stability}) For each irreducible component $M$ as above, the ratio $|M(\BF_q)|/q^{\dim M}$ stabilizes to a limit as the degree of the curve parametrized by $M$ increases.
\end{enumerate}

\smallskip
Over the complex numbers, the first assumption when $C = \BP^1$ is verified in \cite{LT19a} by employing tools from the minimal model program. The second assumption is false in general due to the existence of thin sets, whose geometric characterization is elucidated in \cite{LST22}. Geometric Manin's Conjecture is thence proposed to differentiate the subtlety in the second assumption. In particular, the conjecture singled out a subset of irreducible components of $\Mor(C,X)$ called the Manin components, and the number of Manin components for each curve class are expected to be bounded by a constant. 

\begin{conj}[Geometric Manin's Conjecture (IV)]
    \cite[Section 4]{LRT25} \label{GMC4} Let $X$ be a smooth Fano variety over $\BC$ with Brauer group $Br(X)$. Let $C$ be a smooth connected genus $g$ curve. There exists a curve class $\alpha \in Nef_1(X)_{\BZ}$ such that for all $\beta \in \alpha + Nef_1(X)_{\BZ}$, there are exactly $|Br(X)|$ many Manin components in $\Mor(C, X,\beta)$.
\end{conj} 

Beside the Manin components, the other components could exhibit  pathological behaviors, as they may parametrize non-dominant families of curves, or dominant families such that the evaluation maps have disconnected fibres. In \cite{LRT25}, the authors defined two candidates for pathological components: the accumulating components and the exceptional components. In particular, the accumulating components are geometric in nature, whereas the exceptional components are more suitable for points counting problems over $\BF_q$. In the case of Fano fibrations, \cite{LRT25} completely described the behavior of the pathological families of curves using Fujita's $a$-invariant and showed that the relevant families come from a bounded family of accumulating maps.

The third assumption can be understood as a statement of homological or motivic stability, and less is known compared to the previous two assumptions. We refer the curious readers to \cite[Section 4]{LRT25} for a comprehensive discussion on Batyrev's heuristic and a detailed formulation of Geometric Manin's Conjecture.

\subsection{\label{resandstra} Results and strategies}
To tackle Theorem \ref{IrredM10}, we divide the irreducible components into two categories, the free components and the non-free components. This terminology is based on the notion of a free curve.

\begin{definition}[Free curve] 
    Let $C$ be a smooth projective curve. We say a morphism $f: C\to X$ is a free curve if $f^*T_X$ is generically globally generated and $h^1(C, f^*T_X) = 0$. We say $f$ is non-free if any of the two conditions fails.
\end{definition}

Let $M$ be an irreducible component of either $\Mor(C,X)$ or $\widetilde{M}_{1,0}(X)$. Since being free is an open condition, we say $M$ is a free component if it generically parametrizes free curves. Otherwise, we say $M$ is a non-free component. It turns out that when $X$ is a smooth cubic threefold, a Manin component is the same as a free component, and the pathological components are the non-free components. 

\begin{definition} [$a$-invariant] \label{defa-inv} 
    Let $X$ be a smooth projective variety and let $H$ be a big and nef line bundle on $X$. The $a$-invariant, or Fujita invariant, of the pair $(X,H)$ is defined as
    \[
        a(X, H) = \min\ \{ t\in\BR\ |\ K_X + tH \in \Eff^1(X) \},
    \]
    where $K_X$ is the canonical bundle of $X$ and $\Eff^1(X)$ is the pseudo-effective cone of divisors.
\end{definition}

\smallskip
In general, if $M\subset \Mor(C,X)$ is a non-free component, we expect there to be a generically finite morphism $f:Y\to X$ such that $f$ is not both dominant and birational, $a(Y,-f^*K_X)\geqslant a(X,-K_X)$, and an irreducible component $N \subset \Mor(C,Y)$ such that $f_*$ induces a dominant map $f_*: N \to M$. We can summarize the prediction of Geometric Manin's Conjecture as follows:

\smallskip
\begin{enumerate}[\hspace{5mm}1.]
    \item For sufficiently large $e$, there is a unique free component in $\Mor(C,X)$ parametrizing morphisms of degree $e$.
    \item Non-free components arise from morphisms to $X$ that do not decrease the $a$-invariant.
\end{enumerate}

\noindent We remark that free components enjoy the property that the restricted tangent bundle $f^*T_X$ for a general free curve is sufficiently positive, whereas non-free components parametrize curves whose $f^*T_X$ is unbalanced in the sense of its Harder-Narasimhan filtration.

\subsubsection{\textbf{Non-free components}}

It is straightforward to prove the space of degree $e$ covers of lines forms a non-free component. The crucial step towards showing these are the only non-free components lies in understanding the $a$-invariant of subvarieties of $X$.

\begin{theorem} \label{Classifya-inv} 
    Let $X$ be a smooth cubic threefold and $H$ its hyperplane section. If $S$ is a surface contained in $X$, then 
    \begin{enumerate}
        \item either $a(S, H|_S) \leqslant 1$, or
        \item $a(S, H|_S) = 2$, and a resolution $\tilde{S}$ of $S$ factors through the universal family of lines on $X$, i.e. there exists a morphism $g: \tilde{S}\to \CU$, where $\CU$ is the universal family of lines on $X$.
    \end{enumerate}   
\end{theorem}

Equipped with this theorem, the idea to show the non-existence of other non-free components can be summarized as the following. Suppose $M$ is a non-free component not parametrizing covers of lines. The unbalanced restricted tangent bundle allows us to construct a sublocus $W$ in $M$ such that the universal family above $W$ is mapped to a proper subvariety $S$ in $X$. We then use computations of the $a$-invariant to argue that such $S$ must come from the universal family of lines, contradicting the assumption on $M$. In particular, we obtain the following classification of non-free components in Section \ref{Non-free components}.

\begin{theorem}
    \label{non-free-main} Let $X$ be a smooth cubic threefold. The unique non-free component in $\widetilde{M}_{1,0}(X,e)$ generically parametrizes degree $e$ covers of lines on $X$.
\end{theorem}

\subsubsection{\textbf{Free components}}

As mentioned in the introduction, the key to classify free components is to prove a Movable Bend-and-Break type result, developed in \cite{BLRT22} for smooth Fano threefolds in the case $C = \BP^1$. Roughly speaking, we set up an incidence condition by fixing general points and general curves on $X$. We then look at the locus that satisfies this incidence condition in an irreducible component and examine its limiting behavior on the boundary. We apply this approach to craft a genus one version of Movable Bend-and-Break when $X$ is a smooth cubic threefold.

\begin{named}{\rm \textbf{Theorem} \ref{Thm:MBB} (Movable Bend-and-Break)}
Let $X$ be a smooth cubic threefold. For $e\geqslant 5$, any free component $M\subset\widetilde{M}_{1,0}(X,e)$ contains a stable map $f: Z = Z_1 + Z_2 \to X$ satisfying the following conditions:
\begin{itemize}
    \item $Z_1$ is a smooth genus one curve and $f|_{Z_1}$ is free.
    \item $Z_2$ is a rational curve and $f|_{Z_2}$ is free.
    \item $f|_{Z_i}$'s are birational onto their images.
\end{itemize}
\end{named}

Movable Bend-and-Break can be viewed as a refined version of the celebrated Mori's Bend-and-Break in special settings, as it gives a precise description of the shape of the broken curve. Knowing the broken curve allows us to perform a movement of breaking and deforming. More specifically, by assuming the irreducibility of lower degree free components and using classifications of low degree curves in \cite{HRS05}, the breaking and deforming move allows us to inductively prove the irreducibility of higher degree free components. This strategy is initiated by \cite{HRS04} and is also reminiscent to the story of the connectedness of $\overline{M}_{g,n}$. We perform this movement in Section \ref{Free components}.

There are two essential ingredients inherited in the setup of the Movable Bend-and-Break. First of all, the normal sheaf for a general stable map $f:E\to X$ in a free component is locally free when the degree is sufficiently large. The proof relies on the insight that the slope stability of a semistable torsion-free sheaf $\CF$ on $X$ is ``close" to the slope stability of its restriction to a nef curve class. In particular, this allows us to control the Harder-Narasimhan filtration of $f^*T_X$, which in turn allows us to deduce that $f$ is an immersion, i.e. its normal sheaf is locally free. The technical inputs are referred to as the Grauert-M\"ulich type results, which we briefly recall in Section \ref{SectionGM}.

Secondly, we are able to utilize Theorem \ref{Classifya-inv} to obtain the balanced condition of the normal bundle in Proposition \ref{bal}. The balanced condition is vital as it enables us to fix the maximum number of general points when deforming the stable maps in Proposition \ref{Prop:MBB}, and the number of general points a curve passes through determines its free-ness. 

\smallskip

Finally, we establish the irreducibility of the morphism space in Theorem \ref{IrredMor} by identifying $\Mor(E,X)$ with a sublocus of the fibre of the forgetful morphism $\overline{\CM}_{1,1}(X)\to\overline{M}_{1,1}$. We use Stein factorization to study the irreducibility of a general fibre.

\subsection{\label{comparisongenus0} Comparison to genus zero map.}
We describe some challenges arising from studying the space of higher genus curves. The first prominent difference between rational curves and higher genus curves presents in the definition of free curve. Let $f: C\to X$ be a free curve. When $C = \BP^1$, the condition on global generation of $f^*T_X$ implies the vanishing of $h^1(C, f^*T_X)$. But this is not true when the genus is greater than $0$. 

\begin{example} 
    Let $C$ be a smooth connected curve of genus $g$, and consider a family of plane curves $f:C\to \BP^2$ of degree $d$ dominating $\BP^2$. If $2g > d^2$, then a straightforward computation shows that $h^1(C, f^*T_{\BP^2}) > 0$ for degree $d \geqslant 5$, hence such dominant family cannot be free.
\end{example}

\smallskip
Moreover, the geometric interpretation of free and non-free curve differs. If a family of rational curves maps dominantly onto the variety, then a general member of this family is a free curve. But this implication is not true in the higher genus scenario, as we will see in Subsection \ref{Unifamline}. 

Another subtlety is the more complicated structure of the Kontsevich moduli space of higher genus stable maps. For example, the locus of stable maps with a contracted component of genus at least one in the domain could form irreducible components, and they could meet the main components in low codimension. This prompts us to take further care when performing Movable Bend-and-Break. 

\begin{example}
    \cite{VZ07} Consider the Kontsevich moduli space $\overline{M}_{1,0}(\BP^2, 3)$ of cubic genus one curves on $\BP^2$. There is a unique main component $M_0$ of dimension $9$ generically parametrizing smooth cubics. On the other hand, there are also irreducible components parametrizing maps with reducible domain. For example, the locus parametrizing stable maps of the form $f:E+\BP^1\to \BP^2$ with $f|_E$ constant induces an irreducible component $M_1$ of dimension $10$; and $M_1\cap M_0$ is a codimension $1$ subset in $M_0$.
\end{example}

\subsection{Previous works}
There are several varieties where the components of moduli space of genus one curves have been explored. When the target is a Grassmannian, the morphism space is studied in \cite{Bru87} and the Kontsevich moduli space is studied in \cite{Qua23}. For certain homogeneous spaces, the morphism space is addressed in \cite{Per12} and \cite{PP13}. More generally, there are also works concerning higher genus curves. For example, the Hilbert scheme of higher genus curves on smooth quadric is computed in \cite{Bal89}, the connectedness of Kontsevich moduli space of homogeneous spaces is established in \cite{BK01}, and the morphism space for smooth hypersurfaces of very low degree is analyzed in \cite{Has24}.

\medskip
\noindent\textbf{Leitfaden:} The paper is arranged as follows. In Section \ref{Preliminaries}, we provide basic definitions and tools used throughout the paper. In Section \ref{a-inv}, we collect results on $a$-invariants of surfaces in smooth cubic threefolds and settle Theorem \ref{Classifya-inv}. Section \ref{Non-free components} and \ref{Free components} form the core of the paper. Specifically, in Section \ref{Non-free components}, we prove irreducibility of the components arising from covers of lines. We also show there is no non-free component parametrizing dominant family of curves. In Section \ref{Free components}, we prove the Movable Bend-and-Break Theorem \ref{Thm:MBB}. Combining with irreducibility results on the space of rational curves, we conclude our main Theorem \ref{IrredM10} by inductively proving the irreducibility of the space of free degree $e$ genus one curves.

\medskip
\noindent\textbf{Acknowledgments:}  
The author would like to express his sincere gratitude towards his advisor Brian Lehmann for his suggestion of the topic, his precious advice, his invaluable patience, his tremendous encouragement, and his detailed comments and feedback throughout the project. Special thanks goes to Eric Jovinelly for pointing out an argument leading to the proof of Theorem \ref{IrredMor}, for suggesting a few references, and for his careful reading of the article. The author would also like to thank Benjamin Church, Eric Jovinelly, Fumiya Okamura, and Eric Riedl for discussions, and Ting Gong and Fumiya Okamura for providing useful feedback on an earlier draft of the paper.

\section{\label{Preliminaries} Preliminaries}

We work over $\BC$ throughout the paper. The variety we have in mind will be a separated, integral scheme of finite type over $\BC$. Curves are always assumed to be projective and connected, and we let $E$ denote a smooth genus one curve in the rest of the paper.

Let $X$ be a projective variety. Let $N^1(X)_{\BR}$ denote the space of $\BR$-Cartier divisors on $X$ up to numerical equivalence. We denote by $\Eff^1(X)$ and $\Nef^1(X)$ the pseudo-effective and nef cones of divisors respectively. Dually, let $N_1(X)_{\BR}$ denote the space of real $1$-cycles up to numerical equivalence, and we denote by $\Eff_1(X)$ and $\Nef_1(X)$ the pseudo-effective cone and the nef cone of curves respectively. The intersection pairing $N_1(X)_{\BR} \times N^1(X)_{\BR} \to \BR$ induces isomorphisms $\Eff^1(X) \cong \Nef_1(X)^\vee$ and $\Eff_1(X) \cong \Nef^1(X)^\vee$.

{\definition\label{free} Let $X$ be a smooth projective variety, $C$ be an irreducible projective curve, and $f: C\to X$ be a morphism. Let $p\in C$ be a general point. We say $f$ is a free (resp. very free) curve if $f^*T_X$ (resp. $f^*T_X(-p)$) is generically globally generated and $h^1(C, f^*T_X) = 0$ (resp. $h^1(C, f^*T_X(-p)) = 0$). If any one of the conditions fails, we say $f$ is non-free. }

{\remark In the literature (ex. \cite{Kol96}), a curve $f: C\to X$ is free if $f^*T_X$ is globally generated and $h^1(C, f^*T_X) = 0$. Our notion of free curve is less restrictive and is more useful in classifying the irreducible components later on.}

\subsection{Deformation theory of curves}\hfill

\smallskip
\noindent Let $X$ be a smooth projective variety. The Kontsevich moduli stack of stable maps $\overline{\CM}_{g,0}(X)$ is a Deligne-Mumford stack whose coarse moduli space is a projective scheme $\overline{M}_{g,0}(X)$. We briefly discuss its deformation theory following \cite{BM96} and \cite{BF97}. In Section \ref{Free components}, we will prove more specific results to maintain the flow of exposition.

Let $f: Z \to X$ be a stable map of genus $g$. Assume that $f$ is an immersion on an open neighborhood of each node of $Z$. The deformation theory of $f$ is controlled by the normal sheaf $N_f$ which fits into an exact sequence:
\[
0 \to \CExt^1_{\CO_Z}(Q, \CO_Z) \to N_f \to \CHom_{\CO_Z}(K, \CO_Z) \to 0,
\]
where $K$ and $Q$ are the kernel and cokernel of the complex $f^*\Omega^1_X \to \Omega^1_Z$. The space of first order deformations of $f$ is $H^0(Z, N_f)$ and its obstruction space is $H^1(Z, N_f)$. In particular, $f$ is a smooth point in $\overline{M}_{g,0}(X)$ if $h^1(Z, N_f) = 0$. We will mainly be dealing with the case when $f$ is an immersion, in which case $Q = 0$. 

Assume $f$ is a stable map of genus one. When the dual graph of the domain is a tree, the following lemma tells us when a stable map is a smooth point in the Kontsevich moduli space and when a smoothing of such a stable map exists:

{\lemma \cite[Lemma 2.3]{HRS05} \label{letssmoothNf} Let $f:Z\to X$ be a stable map of genus one. Assume $f$ is an immersion and $Z = Z_1 + Z_2$ is a union of two nodal curves $Z_1$ and $Z_2$. If $p = Z_1\cap Z_2$ and $H^1(Z_1, N_{f|_{Z_1}}(-p)) = H^1(Z_2, N_{f|_{Z_2}}) = 0$, then $[f]$ is a smooth point in $\overline{M}_{1,0}(X)$ and $f$ can be deformed to a stable map $f': E \to X$, where $E$ is a smooth genus one curve.}

\smallskip
Next, we derive several useful results on the evaluation map of the morphism scheme $\Mor(C,X)$ for a smooth projective curve $C$. We refer the reader to \cite[Chapter 2]{Deb01} for notations and \cite{Kol96} for a detailed introduction to this scheme in a more general setting. The following propositions are translations of the cited results for rational curves in \cite[Chapter 4]{Deb01} to the genus one setting, and we omit the proofs.

{\prop \cite[c.f. Proposition 4.9]{Deb01} \label{dominantggg} Let $X$ be a smooth projective variety of dimension $n$. Let $M\subset \Mor(E,X)$ be an irreducible component. Suppose $M$ parametrizes a dominant family of curves on $X$, i.e,
\[
ev: E \times \Mor(E, X) \to X
\]
is dominant. Then for a general member $g:E\to X$ of $M$, the restricted tangent bundle $g^*T_X$ is generically globally generated.
}

\smallskip
Note that if we fix a non-empty finite subscheme $K$ of $E$ such that the evaluation map of $\Mor(E, X, s|_K)$ is dominant, the conclusion of the above theorem does not hold, i.e. $g^*T_X(-K)$ may not be generically globally generated for a general $[g] \in \Mor(E, X, s|_K)$. Nevertheless, we could still obtain a similar result if the decomposition of the bundle is special. We need the following result when $X$ is a smooth projective threefold.

{\prop \label{nondom} Let $X$ be a smooth projective threefold. Let $M \subset \Mor(E,X)$ be an irreducible component. Suppose $M$ sweeps out a dominant family of curves on $X$ and for a general member $s:E\to X$ of $M$, the restricted tangent bundle $s^*T_X = \CO_E\oplus \CO_E\oplus \CL$. Let $p$ be a point on $E$. The following evaluation map is not dominant:
\[
ev: E\times \Mor(E, X, s|p) \to X.
\] 
}

\subsection{Smooth cubic threefolds} \hfill

\smallskip
\noindent A smooth cubic threefold $X$ is a smooth degree $3$ hypersurface in $\BP^4$. This implies $\Pic(X) = \BZ H$ and $-K_X = 2H$. Clemens and Griffiths answered the rationality question of smooth cubic threefolds in \cite{CG72} by introducing the method of Intermediate Jacobian. Along the way, they established detailed descriptions of the Fano surface of lines on smooth cubic threefolds.

{\definition \cite{CG72} Let $X$ be a smooth cubic threefold. The Fano surface of lines $V$ on $X$ is defined to be the Hilbert scheme of lines on $X$. Given a line $L\subset X$, we say
\begin{enumerate}
\item $L$ is a line of first type if $N_{L/X} = \CO_L \oplus \CO_L$, and
\item $L$ is a line of second type if $N_{L/X} = \CO_L(1) \oplus \CO_L(-1)$.
\end{enumerate}}

{\theorem \cite{CG72} The Fano surface of lines $V$ is a smooth surface of general type. Through a general point on $X$, there are $6$ lines of the first type. Lines of the second type form a codimension one locus in $V$.}

\smallskip
The geometry of $V$ is well studied. For example, the Albanese morphism of $V$ is an embedding:

{\prop \cite[Cor. Parag. 4]{Bea81} \label{noratonV} The map $V\rightarrow \text{Alb}(V)$ is an embedding. In particular, $V$ contains no rational curves.}

{\corollary \cite[Lemma 3.8]{GK19} \label{GK19Lemma3.8} Any non-constant morphism $f:E\to V$ from a smooth genus one curve to $X$ has image a smooth genus one curve.}

\smallskip
Moreover, for any smooth cubic threefold, the number of smooth genus one curves on $V$ is finite and at most 30 \cite[Theorem 26]{Rou09}, and $V$ does not contain any smooth genus one curve for a general $X$:

{\prop \cite[Theorem 26]{Rou09} \label{finiteEonV} Let $n_B$ be the number of smooth genus one curves on $V$. Then Fano surfaces with $n_B > 0 $ form a $7$-dimensional family in the $10$-dimensional moduli space of Fano surfaces. In particular, a general Fano surface $V$ contains no smooth genus one curve. }

\smallskip
Let $\CU$ denote the universal family over $V$:
\[\begin{tikzcd}
	{\mathcal{U}} & X \\
	V
	\arrow["\pi"', from=1-1, to=2-1]
	\arrow["ev", from=1-1, to=1-2]
\end{tikzcd}\]
We have a nice description of the tangent bundle $T_V$:

{\theorem[Tangent bundle theorem] \cite[Theorem 3]{Rou11}\cite[Theorem 12.37]{CG72} \label{TangentBundleTheorem}
There is an isomorphism 
\[
\pi_* ev^* \CO_X(1) \cong \Omega_V,
\]
where $\Omega_V$ is the cotangent bundle. In particular, we can identify $H^0(X,\CO_X(1))$ with $H^0(V, \Omega_V)$, the variety $\BP^4$ with $\BP(H^0(V,\Omega_V))$, and the universal family $\CU$ over $V$ is equal to the projectivized tangent bundle $\BP(T_V)$ over $V$. Here, points of the projectivization corresponds to 1-dimensional quotients.}

\smallskip
Under the evaluation map, the lines on $X$ parametrized by smooth genus one curves in $V$ sweep out a cone on $X$. Each smooth genus one curve therefore corresponds to the classical notion of an Eckardt point on $X$.

{\lemma \cite[Proposition 10]{Rou09} \label{self-intersectionB} Suppose there is a smooth genus one curve $B$ on $V$. The self-intersection number of $B$ is $B^2 = -3$. Furthermore, $ev(\pi^{-1}(B))$ is a cone over a plane cubic on $X$. In particular, the cone is the intersection of $X$ with the hyperplane $T_{X,p_B}$ where $p_B$ is the vertex of the cone. }

\smallskip
The Fano surface of lines also provides an inductive way of studying the Hilbert schemes of low degree curves on $X$ via residuation: in \cite{HRS05}, the authors proved the following results with the aim of studying the birational geometry of these Hilbert schemes:

{\theorem \cite{HRS05} \label{basecase} Let $\CH^{d,1}(X) \subset Hilb_{dt}(X)$ denote the open subscheme parametrizing smooth curves of degree $d$ and genus 1. 
\begin{enumerate}[(i)]
    \item The Hilbert scheme of genus one cubic curves $Hilb_{3t}(X)$ on $X$ is isomorphic to the Grassmannian $\BG(2,4)$. In particular, it is smooth and connected of dimension $6$.
    \item $\CH^{4,1}(X)$ is smooth and connected of dimension $8$.
    \item $\CH^{5,1}(X)$ is irreducible of dimension $10$.
\end{enumerate}
}

\smallskip
Since the irreducible main components of the Kontsevich moduli space of genus one curves are birational to their corresponding Hilbert schemes, this result will serve as our base case of induction in the study of free components. Let $H$ be the generator of the Picard group of a smooth cubic threefold. We use $\overline{M}_{g,0}(X,e)$ to denote the union of irreducible components in $\overline{M}_{g,0}(X)$ parametrizing stable maps of $H$-degree equal to $e$.

{\corollary \label{base1} For $e\in\{3,4,5\}$, there is a unique irreducible component in $\overline{M}_{1,0}(X,e)$ generically parametrizing stable maps with irreducible domains that are birational onto the image.
}

{\proof Let $M \subset \overline{M}_{1,0}(X,e)$ be an irreducible component generically parametrizing stable maps with irreducible domains that are birational onto their images. Let $f:E\to X$ be a general member of $M$. The image of $f$ is either a smooth genus one curve or a singular curve of genus greater than one. The latter case is not possible by \cite{HRS05}. Hence $M$ parametrizes stable maps with smooth images. The universal property of Hilbert schemes implies that there is a birational map $\phi: M \dashrightarrow \CH^{e,1}$, where $\CH^{e,1}$ is the Hilbert scheme of degree $e$ genus one curves on $X$. Hence the irreducibility of $M$ follows from the irreducibility of $\CH^{e,1}$.
\qed}

\smallskip
The following theorem gives the irreducibility of the Kontsevich moduli space of rational curves on any smooth cubic threefold, which is crucial in the deformation argument in our main theorem.

{\theorem\cite[Theorem 62]{Sta00}\label{Starr} Let $X$ be a smooth cubic threefold. For each degree $e\geqslant 2$, there are exactly two irreducible components of dimension $2e$ in $\overline{M}_{0,0}(X,e)$, of which one generically parametrizes maps birational onto their images, and the other generically parametrizes degree $e$ covers of lines.}

\subsection{Vector bundles on smooth genus one curves}\hfill

\smallskip
\noindent In this section, we collect and make explicit some results about vector bundles on smooth genus one curves. The following results are direct consequences of \cite{Ati57}. 

{\prop \label{gggindecompvectbundle} Let $\CE$ be an indecomposable vector bundle on a smooth genus one curve $E$ of rank $r$ and degree $d$. Suppose $\CE$ is not the trivial line bundle $\CO_E$. Then $\CE$ is generically globally generated if and only if $d \geqslant r$.}

{\corollary \label{gggh^1} Let $\CE$ be a generically globally generated vector bundle on a smooth genus one curve $E$. Suppose furthermore that $h^1(E, \CE) = k$. Then $\CE$  contains $k$ copies of $\CO_E$ as direct summands.
}

{\theorem \label{gggvectbundle} Let $\CE$ be a vector bundle of degree $d > 0$ and rank $r$ on a smooth genus one curve $E$, then $\CE$ is globally generated (resp. generically globally generated) if and only if the following condition holds:
\begin{enumerate}[(i)]
\item If $\CE$ is indecomposable, then $d > r$ (resp. $d\geqslant r$).
\item If $\CE$ is decomposable, then each summand is globally generated (resp. generically globally generated).
\end{enumerate}
}

\smallskip
We also know when an indecomposable vector bundle is semistable: 

{\theorem \cite{Ati57}\cite[Appendix A]{Tu93} Let $\CE$ be an indecomposable vector bundle of degree $d>0$ and rank $r$ on a smooth genus one curve E. Then $\CE$ is semistable if and only if $d\geqslant r$.
}

\subsection{Grauert-M\"ulich \label{SectionGM}}\hfill

\smallskip
\noindent Given a semistable vector bundle $\CE$, the theorem of Grauert and M\"ulich \cite{GM75} describes the Harder-Narasimhan filtration of the restriction of $\CE$ to lines in $\BP^n$. In \cite{LRT23a} and \cite{LRT25}, the authors extend the result to study non-free sections (resp. curves) of Fano fibrations (resp. Fano varieties). These series of results are refered to as ``Grauert-M\"ulich". 
In this section, we fix some basic definitions and review a key theorem developed in \cite{LRT25}, which we will apply afterwards to study the positivity of the restricted tangent bundles of $X$. We adopt the following notions of slope stability with respect to movable curve classes developed in \cite{CC11}.

{\definition Let $Z$ be a smooth projective variety and $\alpha \in Nef_1(Z)$. Let $\CE$ be a torsion-free sheaf on $Z$. We define the slope of $\CE$ with respect to $\alpha$ as
\[
\mu_\alpha(\CE) = \frac{c_1(\CE)\cdot \alpha}{rk(\CE)}.
\]
We say that $\CE$ is $\alpha$-semistable if for every non-zero torsion free subsheaf $\CF \subset \CE$, we have $\mu_\alpha(\CF) \leqslant \mu_\alpha(\CE)$. If $Z$ is a smooth projective curve, we take $\alpha$ to be the fundamental class of $Z$ and simply write $\mu(\CE)$.
}

{\definition \cite[Definition 2.11.]{LRT25} Let $Z$ be a smooth projective variety and let $\alpha \in Nef_1(Z)$. Suppose $\CE$ is a torsion-free sheaf of rank $r$. Write 
\[
0 = \CF_0 \subset \CF_1 \subset \dots \subset \CF_k = \CE
\]
for the $\alpha$-Harder Narasimhan filtration of $\CE$. The slope panel $SP_\alpha(\CE)$ is the $r$-tuple of rational numbers (in non-increasing order):
\[
SP_\alpha(\CE) = (\underbrace{\mu_\alpha(\CF_1/\CF_0), \dots}_{rk(\CF_1/\CF_0) \ \text{copies}}, \dots, \underbrace{\mu_\alpha(\CF_k/\CF_{k-1}), \dots}_{rk(\CF_k/\CF_{k-1}) \ \text{copies}}).
\]
We denote by $\mu_\alpha^{\max}(\CE)$ the maximal slope of any torsion-free subsheaf, i.e. $\mu_\alpha^{\max}(\CE) = \mu_\alpha(\CF_1)$ and $\mu_\alpha^{\min}(\CE)$ the minimal slope of any torsion-free quotient, i.e. $\mu_\alpha^{\min}(\CE) = \mu_\alpha(\CE/\CF_{k-1})$.
If $Z$ is a smooth projective curve, then we take $\alpha$ to be the fundamental class of $Z$ and simply write $SP(\CE)$.
}

{\theorem[Grauert-M\"ulich] \cite[Theorem 6.5]{LRT25} \label{GM} Let $X$ be a smooth projective variety. Let $W$ be a variety equipped with a generically finite morphism $W \to \overline{M}_{g,0}(X)$. Let $p : U_W \to W$ denote the universal family over $W$ with evaluation map $ev_W : U_W \to X$. Denote by $C$ a general fibre of $p$ with an induced morphism $s: C\to X$. Assume that $C$ is smooth, that $ev_W$ is dominant, that the general fibre of $ev_W$ is irreducible, and that $C$ is contained in the locus where $ev_W$ is flat.

Suppose that $\CE$ is a torsion free sheaf on $X$ of rank $r$ that is semistable with respect to $s_*C$. Write the Harder Narasimhan filtration of $s^*\CE$ as 
\[
0 = \CF_0 \subset \CF_1 \subset \dots \subset \CF_k = s^*\CE.
\]
Let $t$ be the length of the torsion part of $N_s$, let $\CG$ be the subsheaf of $(N_s)_{tf}$ generated by global sections, let $M_{\CG}$ denote the syzygy bundle of $\CG$, and let $V$ be the tangent space to $W$ at $s$. Let $q$ be the dimension of the cokernel of the composition
\[
V \to T_{\overline{M}_{g,0}(X),s} = H^0(C,N_s) \to H^0(C, (N_s)_{tf}).
\]
Then for every index $1\leqslant i \leqslant k-1$, we have
\[
\mu(\CF_i/\CF_{i-1}) - \mu(\CF_{i+1}/\CF_i) \leqslant (q+1)\mu^{max}(({M_\CG})\smvee) + t.
\]
}

\smallskip
In our case, we are interested in bounding the difference between slopes of successive terms in the Harder Narasimhan filtration of $s^*T_X$ over a smooth genus one curve $E$.

\section{\label{a-inv} Fujita's \textnormal{$a$}-invariant}

\noindent In this section, we prove Theorem \ref{Classifya-inv} which describes the surfaces in $X$ with $a$-invariants greater than $1$. Recall from Definition \ref{defa-inv} that given a smooth projective variety $X$ and a big and nef line bundle $L$ on $X$, the $a$-invariant of the pair $(X,L)$ is defined as:
\[
a(X, L) = \min\ \{ t\in\BR\ |\ K_X + tL \in \Eff^1(X) \}.
\] 

\smallskip
Roughly speaking, the $a$-invariant measures how far away the canonical divisor is from the pseudo-effective cone. When $L$ is nef but not big, we set $a(X,L) = \infty$.
The $a$-invariant is also a birational invariant by the next proposition. When $X$ is singular and $L$ is a big and nef Cartier divisor, we define $a(X,L) = a(X', \phi^*L)$, where $\phi:X'\to X$ is any resolution of singularity.

{\prop \cite[Proposition 7]{HTT15} Let $\phi: X' \to X$ be a birational morphism of smooth projective varieties and $L$ a big and nef line bundle on $X$. Then we have $a(X', \phi^*L) = a(X, L)$.}

\smallskip
When $X$ is a uniruled variety, the following theorem allow us to relate the $a$-invariant of $X$ with the $a$-invariant of a subvariety of $X$.

{\definition Let $X$ be a smooth projective variety and $L$ a big and nef line bundle on $X$, we say the pair $(X,L)$ is adjoint rigid if $\kappa(K_X + a(X,L)L) = 0$, where $\kappa(K_X + a(X,L)L)$ is the Iitaka dimension of $K_X + a(X,L)L$.}

\smallskip
The key to the proof of Theorem \ref{Classifya-inv} is the classification of projective varieties with large $a$-invariants.

{\theorem \cite[1.3. Proposition.]{Hor10}\cite{Fuj89} \label{a-inv-cor} Let $Z$ be a smooth projective variety of dimension $n$ over an algebraically closed field and let $H$ be a big and nef divisor on $Z$.
\begin{enumerate}
\item If $a(Z, H) > n$, then $a(Z, H) = n+1$ and the pair $(Z,H)$ is birationally equivalent to $(\BP^n, L)$ where $L$ is the hyperplane section.

\item If $a(Z, H)  = n$ and $(Z, H)$ is adjoint rigid, then $(Z,H)$ is birationally equivalent to $(Q, L)$ where $Q$ is a quadric hypersurface, possibly singular, and $L$ is the hyperplane section on $Q$.

\item If $a(Z, H) = n$ and $(Z, H)$ is not adjoint rigid, then up to birational equivalence, the canonical map $\pi: Z\to C$ realizes $Z$ as a $\BP^{n-1}$-bundle over a curve $C$ and $L = \CO_\pi(1)$.

\item If $n-1 < a(Z,H) < n$, then $(Z,H)$ is birationally equivalent to $(\BP_{\BP^2}(\CO(2)\oplus\CO^{n-2}), L)$ with $L = \CO_{Z/\BP^2}(1)$. In this case, the pair $(Z,H)$ is adjoint rigid and $a(Z,H) = n - \frac{1}{2}$.
\end{enumerate}
}

\smallskip
Thus, to establish Theorem \ref{Classifya-inv}, it suffices to rule out Case (a), (b), and (d) one by one. Case (a) is proven in \cite[Lemma 7.2]{LT19a}, and we borrow similar methods to address the other cases. The basic idea is to replace the surface with a smooth model and run the minimal model program to produce a less complicated variety.

{\prop Let $X$ be a smooth cubic threefold. Let $H$ be the hyperplane section on $X$. There is no adjoint rigid pair $(S,H|_S)$ such that $S$ is a surface in $X$ satisfying $a(S,H|_S) = 3/2$ or $a(S,H|_S) = 2$.
}

{\proof Suppose to the contrary that 
there is an adjoint rigid pair $(S,H|_S)$ satisfying the above hypothesis. Then by Theorem \ref{a-inv-cor}, $S$ is birational to either $\BP^2$ or a quadric surface $Q$. Let $Y$ denote either $\BP^2$ or $Q$. Then there is a rational map $f:Y \dashrightarrow S$. 

Let $\nu: S'\to S$ be the normalization and $\psi: Y \dashrightarrow S'$ be a rational map such that $f = \nu\circ \psi$. Denote by $\phi: W \to S'$ a smooth surface obtained by resolving the indeterminacy of $\psi$. We run the $(K_W + a(S', L)\phi^*L)$-MMP, where $L = \nu^*H|_S$. The first step of the MMP contracts a $(-1)$-curve on $W$. Let $C$ be one such curve, we have $K_W\cdot C = -1$ and $(K_W + a(S', L)\phi^*L)\cdot C < 0$. Since $a(S', L) > 1$, we have $\phi^*L\cdot C = 0$. This implies that $C$ is also contracted by $\phi$ because $L$ is ample. Hence, we obtain a morphism $\phi': W\to W'$ such that $\phi$ factors through $\phi'$. Since $W'$ is smooth, we may replace $W$ by $W'$ to start with.

Repeating this process, we arrive at a minimal model $\widetilde{W}$. By adjoint rigidness, we have $-K_{\widetilde{W}} \equiv a(S', L) \phi^*L$, and $\widetilde{W}$ is a smooth weak del Pezzo surface. Suppose there is a $(-1)$-curve $C$ on $\widetilde{W}$, then $-K_{\widetilde{W}} \equiv a(S', L) \phi^*L$ implies that 
\[
 \phi^*L\cdot C = \frac{1}{a(S', L)}.
\]
Since $a(S', L)$ is either $3/2$ or $2$ and $L$ is ample, we see that $\phi^*L\cdot C$ is not an integer. Hence there is no $(-1)$-curve on $\widetilde{W}$. We divide the discussion based on the possibility of $Y$:

\begin{enumerate}[(i)]
    \item $Y = \BP^2$ and $a(S,H|_S) = 3/2$: 
    
    Since there is no $(-1)$-curve on $\widetilde{W}$, we have $\widetilde{W} = \BP^2$ and $L = \CO_{\BP^2}(2)$. Since the only possible morphism from $\BP^2$ to a normal surface is the identity map, we have $\BP^2 = S'$ and $f$ is a morphism. Let $n$ be an integer such that $S\in |nH|$. Then we have 
    \[
    4 = L^2 = \deg f \cdot H|_S^2 = 3n \cdot \deg f.
    \]
    As both sides of the equation are integers, we obtain a contradiction.

    \item $Y = Q$ and $a(S,H|_S) = 2$:
    
    If $Q$ is a smooth quadric surface, we have $\widetilde{W} = Q$. If $Q$ is singular, then $\phi^*L \cdot C = 0$, where $C$ is the $(-2)$-curve obtained from the blowup of the cone point. Hence the morphism $\widetilde{W} \to Q$ is the blowup of the cone point of $Q$. The same degree computation as above gives a contradiction.
\end{enumerate}
\qed}

\smallskip
We conclude that if a surface $S$ admits a generically finite morphism to a smooth cubic threefold $X$ such that $a(S, f^*H) > 1$, then $S$ factors through the universal family of lines on $X$: 

{\prop \label{a-inv-surface} Let $X$ be a smooth cubic threefold. If $S$ is a surface admitting a generically finite morphism $f: S\to X$ with $a(S, f^*H) > 1$, then $a(S,H) = 2$ and $(S,H)$ is birationally equivalent to a $\BP^1$-bundle over a projective curve $C$. In particular, the map $f$ factors through the universal family of lines.}

{\proof By Theorem \ref{a-inv-cor}, it suffices to prove the last statement. By similar argument as the previous proposition, we may replace $S$ by its minimal model and regard it as a $\BP^1$-bundle $\pi: S\to C$. Denote the fibre of $\pi$ by $F$ and let $L$ be the pullback of $H$ under $f: S \to X$, we have $F \cdot (K_{S} + 2L) = 0$ and $F^2 = 0$. By adjunction, this implies that $F\cdot L = 1$, so the fibres are lines, and $S$ factors through the universal family of lines. 
\qed}

\smallskip
Finally, we show that the largest $a$-invariant less than $2$ of a smooth projective threefold is at most $3/2$.

\begin{theorem}\label{a-inv of 3fold}
    Let $X$ be a normal projective threefold with at most terminal $\BQ$-factorial singularities. Suppose $L$ is a big and nef line bundle satisfying $1 < a(X,L) < 2$. Then $a(X,L)\leqslant 3/2$. 
\end{theorem}

\begin{proof}
    Suppose to the contrary that $a(X,L) > 3/2$. We run a $(K_X + \frac{3}{2}L)$-MMP on the pair $(X,L)$. By \cite[Proposition 5.4]{And13}, we obtain a birational map $\phi: X \to X'$ such that there exists a big and nef line bundle $L'$ on $X'$ satisfying $\phi^*L' = L + E$ for some exceptional divisor $E$. By \cite[Proposition 3.5]{And13}, the pair $(X',L')$ satisfies that $a(X', L') \leqslant \frac{3}{2}$. Hence it suffices to prove that $a(X,L) = a(X', \phi_*L)$. 

    Let $r$ be a rational number satisfying $r\geqslant 1$. We have
    \[
    K_X + rL = \phi^*(K_{X'} + r\phi_*L) + \sum_j a_j E_j,
    \]
    where the $E_j$'s are exceptional and $a_j > 0$. If $K_{X'} + r\phi_*L$ is pseudo-effective, then since $E_j \in \Eff^1(X)$, we have $K_X + rL$ is pseudo-effective, so $a(X,L) \leqslant a(X,\phi_*L)$. On the other hand, since $K_{X'} + r\phi_*L$ is the pushforward of $K_X + rL$ under $\phi$, we have $a(X,L) \geqslant a(X,\phi_*L)$. This concludes the proof.
\end{proof}

\section{\label{Non-free components} Non-free components}

Geometric Manin's Conjecture predicts that non-free components should arise from families of maps that factor through subvarieties with $a$-invariant larger than or equal to the $a$-invariant of $X$. By \cite[Lemma 7.2]{LT19a}, there is no subvariety in $X$ with strictly larger $a$-invariant, and any generically finite dominant morphism $Y\to X$ of degree at least $2$ with $a(Y, f^*H) = a(X, H)$ factors through the universal family of lines. Hence we begin by investigating the universal family of lines. We first fix some notations. Let $\CP$ be a property of a morphism. Let $\CM$ be either $\Mor(E,X)$ or $\overline{M}_{1,0}(X)$. We say an irreducible component $M \subset \CM$ is a $\CP$ component if a general member of $M$ has property $\CP$. For example, we say $M$ is a free component if $M$ generically parametrizes free curves. Moreover, if $C$ is a projective curve and $f:C\to X$ is a morphism, the degree of $f$ always refers to the degree of $f^*H$.

\subsection{Universal family of lines \label{Unifamline}}\hfill

\smallskip
\noindent Let $V$ be the Fano surface of lines on $X$ and denote by $\CU$ the universal family of lines:

\[\begin{tikzcd}
	{\mathcal{U}} & X & \BP^4\\
	V
	\arrow["\pi"', from=1-1, to=2-1]
	\arrow["ev", from=1-1, to=1-2]
        \arrow[hookrightarrow, from=1-2, to=1-3]
\end{tikzcd}\]

We first study the case when the morphisms are multiple covers of lines, i.e. those that factor through the fibres of $\pi$. We note that the morphism space $M\subset \Mor(E,\BP^1)$ parametrizing degree $e > 1$ covers is smooth and connected of dimension $2e$. Indeed, let $\CP$ denote the Poincare bundle on $\Pic^e(E) \times E$ and let $\pi$ be the projection $\Pic^e(E) \times E \to \Pic^e(E)$. Let $\CE = (\pi_*\CP)^{\oplus 2}$. Then $M$ is isomorphic to an open subscheme of the projective bundle $\BP_{\Pic^e(E)}(\CE)$ parametrizing basepoint free pencils.

{\theorem \label{coveroflines} Let $W$ denote the locus in $\Mor(E,X)$ that parametrizes degree $e > 1$ covers of lines. Then $W$ is an irreducible non-free component of dimension $2e+2$. }

{\proof Let $f: E \to X$ be a morphism of degree $e$ in $W$ that factors through $g: E\to L$, where $L$ is a line on $X$. Since there is a map $\phi: W \to V$ sending a morphism $[f]$ to the corresponding line in the Fano surface of lines $V$, the dimension of $W$ at $[f]$ is given by $\dim_{[f]}\Mor(E, X) = \dim_{[g]}\Mor(E, L) + \dim V = 2e + 2$. Since the fibre of $\phi$ is the space of degree $e$ covers of lines, which is irreducible by the discussion above, $W$ is irreducible. Next, we analyze the smoothness of $W$. Denote the inclusion of $L$ into $X$ by $i:L\to X$, we have an exact sequence
\[
0 \to T_{L} \to i^*T_X \to N_{L/X} \to 0.
\]
In particular, since $\mu^{\max}(N_{L/X}) < 2$, the exact sequence splits and $i^*T_X = \CO_{\BP^1}(2) \oplus N_{L/X}$.

Suppose $L$ is a line of the first type, we have
\[
f^*T_X = g^*(\CO_{\BP^1}(2)\oplus \CO_{\BP^1} \oplus \CO_{\BP^1}) = g^*\CO_{\BP^1}(2) \oplus \CO_E \oplus \CO_E
\]
and $h^0(E, f^*T_X) = 2e+2$ and $h^1(E, f^*T_X) = 2$. If $L$ is of second type, we have
\[
f^*T_X = g^*(\CO_{\BP^1}(2)\oplus \CO_{\BP^1}(1) \oplus \CO_{\BP^1}(-1)) 
\]
and $h^0(E, f^*T_X) = 3e$ and $h^1(E, f^*T_X) = e$.

In particular, the tangent space has dimension $2e+2$ at lines of first type and dimension $3e$ at lines of second type. Since $\dim W = 2e+2$, $W$ forms an irreducible component in $\Mor(E,X)$. 
\qed}

\smallskip
Varying the moduli of the smooth genus one curve $E$, there is a unique component in the Kontsevich moduli space parametrizing degree $e$ covers of lines:

{\theorem \label{irred of N_e} For $e\geqslant 2$, there is a unique irreducible component $N_e \subset \widetilde{M}_{1,0}(X,e)$ of dimension $2e+2$ parametrizing degree $e$ covers of lines.}

{\proof Let $E$ be a smooth genus one curve and $M_E\subset \Mor(E,X)$ the unique irreducible component parametrizing degree $e$ covers of lines. Let $M$ denote the union of all $M_E$'s as $E$ varies in the moduli space of smooth genus one curves. It suffices to show that the image of the morphism $\phi: M \to \widetilde{M}_{1,0}(X,e)$ is contained in the same irreducible component. Let $f:E\to X$ be a member of $M$ that is a degree $e$ cover of a line of first type. We show that by considering $f$ as a stable map in $\widetilde{M}_{1,0}(X,e)$ under $\phi$, it has a deformation that deforms the moduli of the domain of $f$.

Let $T_{\widetilde{M}, [f]}$ denote the space of first order deformations of $[f]$ in $\widetilde{M}_{1,0}(X,e)$. We have the following exact sequence relating $T_{\widetilde{M}, [f]}$ and deformations of $E$:
\[
\dots \to T_{\widetilde{M}, [f]} \to H^1(E, T_E) \to H^1(E,f^*T_X) \to \dots
\]
Since $f^*T_X = \CO_E\oplus \CO_E \oplus \CL$ for some globally generated line bundle $\CL$, and the proof of Theorem \ref{coveroflines} implies that the map $T_E \to f^*T_X$ sends $T_E$ to the summand $\CL$, we see that $H^1(E, T_E) \to H^1(E,f^*T_X)$ is the zero map. Hence $T_{\widetilde{M}, [f]} \to H^1(E, T_E)$ is a surjection. Since $H^1(E, T_E)$ is the space of deformations of $E$, there is a unique irreducible component $N_e \subset \widetilde{M}_{1,0}(X,e)$ containing the image of $\phi$. Since $\phi$ has relative dimension one, $N_e$ has dimension $2e+2$ by Theorem \ref{coveroflines}.
\qed}

\smallskip
Next, we address the situation when a family of genus one curves factoring through $\CU$ is not contracted by $\pi$. Since there is no rational curve on $V$ by Proposition \ref{noratonV}, we consider the situation when $V$ contains genus one curves. By Corollary \ref{GK19Lemma3.8}, $V$ does not contain any singular curve whose normalization is a smooth genus one curve, and by Proposition \ref{finiteEonV}, there can be only finitely many smooth genus one curves on $V$.

{\prop Suppose there is a smooth genus one curve $B$ on $V$. The restricted tangent bundle $T_V|_B \cong \CO_B \oplus N_{B/V}$.}

{\proof Consider the relative tangent exact sequence:
\[
0 \to T_B \to T_{V}|_B \to N_{B/V} \to 0.
\]
The set of extensions of $N_{B/V}$ by $T_B$ is given by the Ext group
\[
\Ext(T_B, N_{B/V}) = H^1(T_B \otimes N_{B/V}^\vee).
\]
Since $T_B = \CO_B$ and $N_{B/V}$ has degree $-3$ by Lemma \ref{self-intersectionB}, $h^1(T_B \otimes N_{B/V}^\vee) = 0$. Hence the only extension is the trivial extension, so $T_V|_B = \CO_B \oplus N_{B/V}$.
\qed}

\smallskip
Let $B$ be a smooth genus one curve on $V$. The previous proposition shows that the restricted tangent bundle is $T_V|_B \cong \CO_B\oplus N_{B/V}$. By Theorem \ref{TangentBundleTheorem}, we denote by $\CY = \BP(T_V|_B)$ the universal family of lines over $B$:
\[\begin{tikzcd}
	\CY = \BP(T_V|_B) & {\mathcal{U}} & X \\
	B & V
	\arrow["\pi_B"', from=1-1, to=2-1]
	\arrow["i",hook, from=1-1, to=1-2]
	\arrow["\pi"', from=1-2, to=2-2]
	\arrow["i_B",hook, from=2-1, to=2-2]
	\arrow["ev", from=	1-2, to=1-3]
\end{tikzcd}\]

\smallskip
Let $f:E\to \CY$ be a morphism. Denote by $g = \pi_B\circ f:E\to B$ the induced morphism. We can construct the following Cartesian diagram
\[\begin{tikzcd}
	{\CZ = \mathcal{Y}\times_BE} && {\mathcal{Y}} \\
	E && B
	\arrow["{\pi_E}"', from=1-1, to=2-1]
	\arrow["h", from=1-1, to=1-3]
	\arrow["g"', hook, from=2-1, to=2-3]
	\arrow["{\pi_B}",from=1-3, to=2-3]
\end{tikzcd}\]
such that $f$ corresponds to a unique section $\sigma:E\to \CZ$ of the pullback bundle $\pi_E: \CZ \to E$. In particular, we have $\CZ \cong \BP(g^*T_V|_B)$. 

Suppose the morphism $g$ is finite of degree $d$. The pseudo-effective cone of $\CZ$ can be described as follows. Let $C_0$ denote the rigid section  with $C_0^2 = -3d$ and $F$ denote the fibre of $\pi_E$. The pseudo-effective cone of $\CZ$ is 
\[
\Eff^1(\CZ) = \BR_{\geqslant 0}[C_0] + \BR_{\geqslant 0}[F].
\]
In particular, the irreducible curves are represented by either $C_0$ or $aC_0 + bF$ with $a\geqslant 0$ and $b\geqslant 3ad$. The canonical divisor is $K_\CZ = -2C_0 - 3dF$.

{\theorem \label{facunifamline} Let $M \subset \Mor(E, \CU)$ be an irreducible component parametrizing morphisms whose images are not contracted by the projection $\pi:\CU \to V$. Then $M$ does not induce an irreducible component in $\Mor(E, X)$, i.e. the image $ev_*(M)$ is not an irreducible component of $ \Mor(E,X)$.}

{\proof By the discussion above, we may assume $V$ contains a smooth genus one curve $B$ such that every morphism parametrized by $M$ lies above $B$. Let $\CY$ and $\CZ$ be constructed as above. Let $f: E\to \CY$ be a morphism and let $g = \pi_B\circ f: E\to B$. It suffices to consider the situation when $g$ is not a constant map. Let the degree of $g$ be $d$. Let $\sigma:E\to \CZ$ denote the section corresponding to $f$. We have $g^*T_V|_B = \CO_E\oplus g^*N_{B/V}$. A section of $\pi_E$ corresponds to a surjection of $g^*T_V|_B$ onto a line bundle $\CL$:
\[
g^*T_V|_B \to \CL \to 0.
\]
Denote by $S_\CL$ the space of sections of $\pi_E$ corresponding to all surjections $g^*T_V|_B \to \CL$. The dimension of $S_\CL$ is given by
\[
\dim S_\CL = \dim H^0(E, \Hom(g^*T_V|_B, \CL)) - 1 = h^0(E, \CL) + h^0(E, \CL\otimes g^*N_{B/V}^\vee) - 1.
\]

A section of $\CZ$ has the form $C_0$ or $C_0 + bF$ for $b\geqslant 3d$. 
\begin{enumerate}[\hspace{3mm} (i)]
    \item If $\sigma(E) = C_0$, then the image of $E$ in $X$ is the vertex of the cone $Y = ev(\CY)$ by Lemma \ref{self-intersectionB}. This implies that $\dim M = 1$.

    \item If $\sigma(E) = C_0 + bF$ with $b\geqslant 3d$, then $\CL$ is a line bundle of degree $b - 3d \geqslant 0$. We have 
    \[
    \dim S_\CL = h^0(E, \CL) + h^0(E, \CL\otimes g^*N_{B/V}^\vee) - 1 \leqslant b - 3d + 1 + b - 1 = 2b - 3d.
    \]
\end{enumerate}
Next, we compute the actual dimension of $M$. Since the morphisms parametrized by $M$ factors through $\CY$ by definition, we may regard $M$ as a locus inside $\Mor(E,\CY)$. Consider the following morphism
\[
\phi: \Mor(E,\CY) \to \Mor(E,B)
\]
sending $[f]$ to $[\pi_B\circ f]$. The fibre of $\phi$ is the space of sections $S_\CL$, and $\dim \Mor(E,B) = 1$ because any morphism $g:E\to B$ is a composition of translation and isogeny by \cite[Chapter III]{Sil86}. This gives $\dim S_\CL + \dim \Mor(E,B) \leqslant 2d - 3b + 1$. 

Denote by $f'$ the composition $ev\circ f$. We compute the expected dimension of the component in $\Mor(E,X)$ parametrizing $f'$.  Since the image of $C_0$ is the cone point of $Y$ and $F$ is a line on the ruling, we have
\[
\dim_{[f']}\Mor(E,X) \geqslant -K_X\cdot f'_*(E) = 2b.
\]
If $b > 0$, then since $2b > 2b - 3d + 1 \geqslant \dim S_\CL + \dim \Mor(E,B)$, $M$ does not induce an irreducible component of $\Mor(E,X)$. If $b = 0$, then $\deg f' = 0$ and $\dim M = 1$ by (i). Since the space of degree $0$ maps is isomorphic to $\text{Aut}(E)\times X$, which is $4$-dimensional, we see that $M$ again does not induce an irreducible component. 
\qed}

\smallskip Combining Theorem \ref{coveroflines} and Theorem \ref{facunifamline}, we immediately obtain the following:

{\corollary \label{facunifamlinecor} Let $M\subset \Mor(E,X)$ be an irreducible component parametrizing maps factoring through the universal family of lines. Then $M$ is a non-free component parametrizing covers of lines. \qed}

\subsection{Double cover of conics}\hfill

\smallskip
\noindent Suppose $W$ is a locus in $\Mor(E,X)$ generically parametrizing maps that are coverings of lower degree curves on $X$. We would like to understand whether such a locus $W$ forms an irreducible component or a proper sublocus inside some larger components.

{\prop \label{component} Let $M$ be an irreducible component of $\Mor(E,X)$ parametrizing a dominant family of curves with degree $e > 4$. Then $M$ either generically parametrizes covers of lines or maps birational onto their images.}

{\proof Let $W\subset \Mor(E,X)$ be an irreducible locus parametrizing morphisms of degree $e > 4$ that are multiple covers of lower degree curves. Assume $W$ does not parametrize covers of lines. We prove that $W$ has dimension less than the expected dimension $2e$.

Suppose first that $W$ generically parametrizes covers of rational curves. Then a general map $s:E\to X$ factors as a composition of a cover $f:E\to \BP^1$ and a map $g: \BP^1 \to X$ that is birational onto its image. Suppose the degree of $g$ is $e_0$ and the degree of $f:E\to \BP^1$ is $e_1$. Since $e > 4$ and $e_0 > 1$, we have $e_0 + e_1 < e_0 \cdot e_1 = e$. By the Hurwitz formula, we have 
\[
\dim_{[s]} \Mor(E,X) = 2e_0 + 2e_1  < 2e.
\]
Next, suppose $M$ generically parametrizes covers of genus one curves. Then a general map $s:E\to X$ factors as a composition of a cover $f:E\to E'$ and a map $g: E' \to X$ that is birational onto its image. Suppose the degree of $g$ is $e_0$ and the degree of $f:E\to E'$ is $e_1$. The dimension of $\Mor(E,X)$ at $[s]$ can be computed as the sum of the dimensions of  $\Mor(E,E')$ and $\Mor(E',X)$ minus one for automorphisms:
\begin{align*}
\dim_{[s]} \Mor(E,X) &= \dim_{[f]} \Mor(E, E') + \dim_{[g]} \Mor(E', X) - 1\\
                    &\leqslant h^0(E, f^*T_E') + h^0(E', g^*T_X) - 1.
\end{align*}
Since $W$ parametrizes a dominant family of curves, we have $g^*T_X$ is generically globally generated. By Theorem \ref{gggvectbundle}, we have $h^0(E', g^*T_X) \leqslant -K_X\cdot g_*E' + h^1(E', g^*T_X) \leqslant 2e_0 + 3$. Since there is no genus one curve of degree less than $3$ on $X$, we have $e_0 \geqslant 3$ and therefore $2e_0 + 3 < 2e$. Hence we have
\[
\dim_{[s]} \Mor(E,X) \leqslant 1 + 2e_0 + 3 - 1 = 2e_0 + 3 < 2e.
\]
\qed}

\smallskip
On the other hand, when $e = 4$, the locus $M\subset \Mor(E,X)$ parametrizing double covers of conics forms an irreducible free component. We include it for completion.

{\theorem \label{irred of double covers of conics} Let $M\subset\Mor(E,X)$ be the locus parametrizing double cover of conics. Then $M$ forms an irreducible component. In particular, there exists a unique irreducible component of $\widetilde{M}_{1,0}(X,e)$ parametrizing double covers of conics.
}

{\proof Let $i: C\hookrightarrow X$ be a general conic on $X$ and $g:E\to C$ be a double cover. Let $f: E\to X$ be an element of $M$ such that $f = i\circ g$. The dimension of $M$ at $[f]$ is given by
\[
\dim_{[f]} M = \dim_{[g]} \Mor(E, C) + \dim \CH^{2,0},
\]
where $\CH^{2,0}$ denotes the Hilbert scheme of smooth conics on $X$. In particular, $\CH^{2,0}$ is irreducible of dimension 4 by \cite{HRS05}. Since $g$ has degree 2, Hurwitz formula gives $\dim_{[g]} \Mor(E, C) = 4$, so $\dim_{[f]}  M = 8$. Furthermore, by \cite[Lemma 3.2]{HRS05}, the normal bundle of a general conic $C$ is $N_{C/X} = \CO(1)\oplus \CO(1)$. Similar computation as in Theorem \ref{coveroflines} shows that $h^1(E, f^*T_X) = 0$ and the tangent space of $M$ at $[f]$ has dimension $h^0(E, f^*T_X) = 8$, so $[f]$ is a smooth point of $M$. This implies that $M$ forms an irreducible free component. Since the fibres of the morphism $\phi: M \to \CH^{2,0}$ are irreducible, we have that $M$ is irreducible. Finally, the last statement follows from a similar argument in Theorem \ref{irred of N_e}.
\qed}

\subsection{Non-free non-dominant case} \hfill

\smallskip
\noindent In this subsection, we show that there is no non-free component parametrizing a non-dominant family of curves on $X$. We will assume the contrary and derive a contradiction using the $a$-invariant.

{\prop \label{dimension bound for family on surface} Let $Y$ be a smooth projective surface of dimension $2$ and suppose $M$ is an irreducible component of $\Mor(E,Y)$ parametrizing a dominant family of curves. Suppose $\dim M > 2$. Let $s:E\to Y$ denote a general map parametrized by $M$. Then we have
\[
-K_Y \cdot s_*E  \leqslant \dim(M) \leqslant -K_Y\cdot s_*E  + 1.
\]}

{\proof The dimension of $M$ is bounded by $h^0(E, s^*T_Y) = -K_Y\cdot s_*E + h^1(E, s^*T_Y)$, so it suffices to show $h^1(E, s^*T_Y) \leqslant 1$ when $M$ is a non-free component. 

Suppose $M$ is a non-free component parametrizing a dominant family of curves. By Definition \ref{free} and Proposition \ref{dominantggg}, $s^*T_Y$ is generically globally generated with $h^1(E, s^*T_Y) \neq 0$. By Corollary \ref{gggh^1}, $s^*T_Y$ is decomposable containing trivial factors, hence $h^1(E, s^*T_X) \leqslant 2$. Moreover, if $h^1(E, s^*T_X) = 2$, then $s^*T_Y = \CO_E \oplus \CO_E$, but this implies $\dim M \leqslant h^0(E, s^*T_Y) = 2$, contradicting the assumption. Hence, we have $h^1(E, s^*T_Y) \leqslant 1$.
\qed}

{\theorem \label{no non-free non-dominant comp} There is no non-free component of $\Mor(E,X)$ parametrizing a non-dominant family of curves on $X$.}

{\proof Suppose to the contrary that $M\subset \Mor(E,X)$ is a non-free component parametrizing a non-dominant family of curves on $X$. Let $Y$ denote a resolution of the proper subvariety swept out by the universal family of $M$ and let $f$ be the map from $Y$ to $X$. It is enough to consider $Y$ being a surface. Let $N$ denote an irreducible component in $\Mor(E,Y)$ such that $f_*: N\to M$ is dominant. Since covers of lines sweep out a dominant family of curves on $X$, we may assume $M$ does not parametrize covers of lines, in which case the expected dimension of $M$ and $N$ are at least $6$. Hence Proposition \ref{dimension bound for family on surface} implies that the dimension of $N$ is bounded by
\[
-K_Y \cdot C \leqslant \dim(N) \leqslant -K_Y\cdot C  + 1,
\] 
where $C = g_*E$ for a general map $g:E\to Y$. Let $s = f\circ g$. Since $\dim N \geqslant \dim M \geqslant  -K_X\cdot s_*E$, we have $1 \geqslant (K_Y - f^*K_X) \cdot C$, and since $-f^*K_X = 2f^*H$, we can rewrite this inequality in the following way:
\begin{align*}
1 &\geqslant (K_Y +2f^*H)\cdot C \\
	&= (K_Y + (\frac{3}{2} + \epsilon) f^*H + (\frac{1}{2} - \epsilon)f^*H)\cdot C,
\end{align*}
where $\epsilon > 0$ is some constant. We have
\[
1-(\frac{1}{2} - \epsilon)f^*H\cdot C \geqslant (K_Y + (\frac{3}{2} + \epsilon) f^*H)\cdot C.
\]
Since $f^*H\cdot C \geqslant 3$,  we may choose $\epsilon$ small enough such that the left hand side of the inequality is less than $0$. Since the movable cone of curves is dual to the pseudo-effective cone of divisors by \cite[Theorem 0.2]{BDPP13}, this shows that $K_Y + (\frac{3}{2} + \epsilon)f^*H$ is not big, i.e. $a(Y, f^*H) \geqslant \frac{3}{2} + \epsilon$. Theorem \ref{Classifya-inv} implies that $a(Y, f^*H) \geqslant 2$, in which case the map $f:Y\to X$ factors through the universal family of lines on $X$, a contradiction to Corollary \ref{facunifamlinecor}.
\qed}

\subsection{Non-free dominant case} \hfill

\smallskip
\noindent Let $M\subset \Mor(E,X)$ be an irreducible component generically parametrizing non-free curves that are birational onto their images such that the universal family dominates $X$. Let $s:E\to X$ be a general member of $M$. By Proposition \ref{dominantggg}, we have $s^*T_X$ is generically globally generated. We divide our discussion based on the following possible indecomposable summand decompositions of the restricted tangent bundle: 

\begin{enumerate}[\hspace{2mm} (1)]
    \item $s^*T_X = \CO_E \oplus \CF$ and $0 \leqslant \mu^{\min}((N_s)_{tf}) \leqslant 1$. 

    \item $s^*T_X = \CO_E \oplus \CF$ and $\mu^{\min}((N_s)_{tf}) \geqslant 2$.

    \item $s^*T_X = \CO_E \oplus\CO_E\oplus\CL$. 
\end{enumerate}

{\definition We say a non-free component $M \subset \Mor(E, X)$ is of \textit{type i} if a general map of $M$ is birational onto its image and satisfies condition $(i)$ above.}

\smallskip
Note that a map $s:E\to X$ that is birational onto its image has degree at least $3$. A global section computation shows that if $M$ is of type 1 or type 2, then $h^1(E, s^*T_X) = 1$, and if $M$ is of type 3, then $h^1(E, s^*T_X) = 2$. We begin by studying type 1 and type 2 components. In particular, we will apply Grauert-M\"ulich (Theorem \ref{GM}) to treat the type 2 components. This requires us to understand better the geometry of $M$.

{\theorem \label{nfconn} Let $M$ be a type 1 or type 2 component parametrizing a dominant family of curves of degree at least $5$. Then a general fibre of the evaluation map $ev: E\times M \to X$ is irreducible.}

{\proof Suppose to the contrary that the evaluation map has reducible fibre. Let $\widetilde{\CU}_M$ denote the normalization of $E\times M$. Denote by $Y'\to X$ the finite part of the Stein factorization of the map $\widetilde{\CU}_M \to X$ and let $Y$ be a resolution of singularities of $Y'$ equipped with a morphism $f: Y \to X$. Let $N \subset \Mor(E,Y)$ be an irreducible component induced by $M$. Let $s': E\to Y$ be a general map of $N$ and let $s = f\circ s'$. 

Since $M$ parametrizes a dominant family of curves of degree at least $5$ on $X$, we have the inequality:
\[
 -K_X\cdot s_*E \leqslant \dim M = \dim N \leqslant -K_Y\cdot s'_*E + h^1(E, s'^*T_Y) \leqslant -K_Y\cdot s'_*E + 2.
\]
Hence we have $(K_Y - f^*K_X) \cdot s'_*E \leqslant 2$. By Theorem \ref{a-inv of 3fold}, the largest possible $a$-invariant less than $2$ of a smooth projective threefold is $3/2$. Since $-K_X = 2H$, we can rewrite the inequality in the following way
\[
\left(K_Y + (\frac{3}{2} + \epsilon) f^*H + (\frac{1}{2} - \epsilon)f^*H\right)\cdot s'_*E \leqslant 2,
\]
where $\epsilon > 0$ is a sufficiently small constant. We have
\[
\left(K_Y + (\frac{3}{2} + \epsilon) f^*H\right)\cdot s'_*E \leqslant 2 - (\frac{1}{2} - \epsilon)f^*H\cdot s'_*E.
\]
If $f^*H\cdot s'_*E = H \cdot s_*E \geqslant 5$, then the right hand side of the inequality is less than $0$. This implies that $a(Y, f^*H) \geqslant \frac{3}{2} + \epsilon$, and Theorem \ref{a-inv of 3fold} in turn implies that $a(Y, f^*H) \geqslant 2$. Hence $s:E\to X$ factors through the universal family of lines by Theorem \ref{Classifya-inv}. This contradicts Corollary \ref{facunifamlinecor}.
\qed}

\smallskip
We will also need to bound the size of the torsion part of the normal sheaf $N_s$.

{\theorem \cite[Corollary 6.11]{AC81} \label{AC81} Let $X$ be a smooth projective variety and let $C$ be a smooth projective curve. Suppose that $M \subset \Mor(C,X)$ is an irreducible component such that the maps parametrized by $M$ dominate $X$ and the general map parametrized by $M$ is birational onto its image. Let $W = M_{red}$ and $[s]$ be a general member of $W$. Then the image of the map $T_{W,[s]} \to H^0(C, N_s)$ has vanishing intersection with $(N_s)_{tors}$.
}

{\prop \label{N_s-tor} Let $M$ be a type 1 or type 2 component and $s: E \to X$ be a general member of $M$.
\begin{itemize}
    \item If $M$ is of type 1, then the length of the torsion of $N_s$ is at most $2$. 
    \item If $M$ is of type 2, then the length of the torsion of $N_s$ is at most $1$. 
\end{itemize}
}

{\proof
Consider the short exact sequence
\[
0 \to T_E \to s^*T_X \to N_s \to 0.
\]
Suppose that $M$ is a type 1 component. The cokernel of the map $H^0(E, s^*T_X)\to H^0(E, N_s)$ has dimension at most $h^1(E, T_E) = 1$. Since $h^1(E, s^*T_X) = 1$, the dimension of the tangent space $T_{W,[s]}$ has codimension at most $1$ in $H^0(E, s^*T_X)$. By Theorem \ref{AC81}, the image of the map $T_{W,[s]} \to H^0(E, N_s)$ has vanishing intersection with $(N_s)_{tors}$, hence $(N_s)_{tors}$ has length at most $2$.

If $M$ is a type 2 component, then $\mu^{\min}((N_s)_{tf}) \geqslant 2$. Hence $h^1(E, N_s) = 0$ and the map $H^1(E, T_E) \to H^1(E, s^*T_X)$ is a bijection. This implies that $H^0(E, s^*T_X)\to H^0(E, N_s)$ is surjective. The same argument as above shows that $(N_s)_{tors}$ has length at most $1$.
}

\smallskip
The following proposition describes the deformation of maps through fixed general points.

{\prop \label{twisting the strata of kontsevich space does not dominate} Let $M\subset \Mor(E,X)$ be a type 1 or type 2 component. Suppose a general member $[s]$ of $M$ satisfies $(N_s)_{tf} = L_1 \oplus L_2$ for $\deg L_1 \leqslant \deg L_2$. Let $n = 1$ if $M$ is of type 1 and $n=\deg L_1$ if $M$ is of type 2. Then the universal family of the locus parametrizing morphisms through $n$ general points does not dominate $X$.
}

{\proof We may replace $M$ by its image $W$ in $\overline{M}_{1,0}(X)$. Let $W^{\circ}\subset W$ be the locus parametrizing morphisms satisfying the above normal sheaf decomposition. Let $[s]$ be a member of $W^\circ$. Let $p_1,\dots,p_n$ be general points in the image of $s$. Let $V \subset W^{\circ}$ be the locus parametrizing deformations of $s$ through the $p_i$'s. Denote by $U$ the universal family of $V$. 

We show that the evaluation map of $U$ does not dominate $X$. Let $p$ be a general point of $E$ not among the $p_i$'s. Consider the commutative diagram with exact rows:
\[\begin{tikzcd}
	0 & {T_{E,p}} & {(T_U)_{([s],p)}} & {(T_{V,[s]})_{p}} & 0 \\
	0 & {T_{E,p}} & {T_{X,s(p)}} & {(N_s)_p} & 0
	\arrow[from=1-1, to=1-2]
	\arrow[from=1-2, to=1-3]
	\arrow[from=1-2, to=2-2]
	\arrow["{d\pi}", from=1-3, to=1-4]
	\arrow["{d(ev)}", from=1-3, to=2-3]
	\arrow[from=1-4, to=1-5]
	\arrow[from=1-4, to=2-4]
	\arrow[from=2-1, to=2-2]
	\arrow[from=2-2, to=2-3]
	\arrow[from=2-3, to=2-4]
	\arrow[from=2-4, to=2-5]
\end{tikzcd}\]
Since $H^0(E, N_s(-p_1-\dots-p_n))$ is the tangent space at $[s]$ of the locus in $\overline{M}_{1,0}(X)$ parametrizing deformations of $[s]$ through the $p_i$'s, we have $T_{V,[s]} \subset H^0(E, N_s(-p_1-\dots-p_n))$. Hence the right most vertical map has rank at most $1$ because $H^0(E, N_s(-p_1-\dots-p_n)) \to (N_s)_p$ is not surjective. Thus the map $d(ev)$ has rank at most $2$. We conclude that the evaluation map of $U$ does not dominate $X$.
\qed}

\smallskip
We can relate the $a$-invariant of a surface swept out by a non-dominant family of curves to the dimension of this family and the degrees of the curves it parametrizes:

{\prop \label{bounding the a-inv of a surface} Let $W \subset \Mor(E, X)$ be an irreducible locus parametrizing curves of degree $e$. Suppose the family of curves parametrized by $W$ sweeps out a surface $S$ in $X$. If $\dim W \geqslant 3$, then we have 
\[
a(S,H|_S) \geqslant \frac{\dim W - 1}{e}.
\]
}

{\proof Let $S$ denote a resolution of the surface swept out by $W$. Denote by $W_S \subset \Mor(E,S)$ an irreducible component containing the induced image of $W$. Then $\dim W_S \geqslant \dim W$.
Let $L$ be the pullback of the hyperplane section to $S$ and let $C$ denote $s_*E$, where $s:E \to S$ is a general map of $W_S$. Since $\dim W \geqslant 3$, we have
\[
\dim W_S \leqslant -K_S\cdot C + 1
\]
by Proposition \ref{dimension bound for family on surface}. Consider the quantity $(K_S + a(S, L)L)\cdot C$. By \cite[Theorem 0.2]{BDPP13} and our dimension bounds, we have
{\allowdisplaybreaks\begin{align*}
0\leqslant \left(K_S + a(S, L)L\right)\cdot C &= K_S\cdot C + a(S,L)L \cdot C  \\
    & \leqslant -\dim W_S + 1 + a(S,L)\cdot e \\
    & \leqslant -\dim W + 1 + a(S,L)\cdot e.
\end{align*}}
This gives
\[
a(S,L) \geqslant \frac{\dim W - 1}{e}.
\]
\qed}

\smallskip
We are ready to rule out the type 1 components:

{\theorem \label{no type 1} There is no type 1 component parametrizing maps of degree $e \geqslant 5$. }

{\proof Assume to the contrary that there exists a type 1 component $M$. Let $s:E\to X$ be a general member of $M$. Let $q \in X$ be a general point and let $p\in E$ be a point mapping to $q$. Let $W \subset M$ denote an irreducible component of the locus containing $[s]$ and parametrizing maps passing through $q$. By assumption, the normal sheaf $N_s(-p)$ is not generically globally generated, hence Proposition \ref{twisting the strata of kontsevich space does not dominate} implies that the family of curves parametrized by $W$ does not dominate $X$. Since $M$ is a type 1 component parametrizing maps of degree $e\geqslant 5$, we have $h^0(E, s^*T_X(-p)) = 2e - 3$ and $h^1(E, s^*T_X(-p)) = 1$. This gives
\[
\dim W \geqslant h^0(E, s^*T_X(-p)) - h^1(E, s^*T_X(-p)) + 1 = 2e - 3 - 1 + 1 = 2e - 3,
\]
where the plus one in the first inequality accounts for varying the choice of $p$. Then Proposition \ref{bounding the a-inv of a surface} gives
\[
a(S,H|_S) \geqslant \frac{2e-4}{e} = 2 - \frac{4}{e}.
\]
Since $e\geqslant 5$ by assumption, we have $a(S,H|_S) > 1$. Since the choice of $q$ is general, Theorem \ref{Classifya-inv} and Corollary \ref{facunifamlinecor} implies that $W$ parametrizes covers of lines. A contradiction.
\qed}

\smallskip
Next, we focus on the type 2 components. We first establish the necessary ingredients needed in Grauert-M\"ulich.

{\prop \label{evflat} 
Let $M$ be a type 2 component and let $W \subset \widetilde{M}_{1,0}(X)$ be the variety induced by $M$. Then the general fibre of the universal family $\pi: \CU_W\to W$ is contained in the smooth locus of the evaluation map $ev: \CU_W\to X$.
}

{\proof Let $s: E\to X$ be a general map parametrized by $W$ and let $p\in E$ be any point. We have the following commutative diagram with exact rows:
\[\begin{tikzcd}
	0 & {T_{E,p}} & {(T_{\CU_W})_{([s],p)}} & {(T_{W,[s]})_{p}} & 0 \\
	0 & {T_{E,p}} & {T_{X,s(p)}} & {(N_s)_p} & 0
	\arrow[from=1-1, to=1-2]
	\arrow[from=1-2, to=1-3]
	\arrow[from=1-2, to=2-2]
	\arrow["{d\pi}", from=1-3, to=1-4]
	\arrow["{d(ev)}", from=1-3, to=2-3]
	\arrow[from=1-4, to=1-5]
	\arrow[from=1-4, to=2-4]
	\arrow[from=2-1, to=2-2]
	\arrow[from=2-2, to=2-3]
	\arrow[from=2-3, to=2-4]
	\arrow[from=2-4, to=2-5]
\end{tikzcd}\]
Since $M$ is a type 2 component, $N_s$ is globally generated, the last vertical map is surjective. Hence $d(ev)$ is also surjective. Since $[s]$ is general, $\CU_W$ is smooth along the fibre $\pi^{-1}([s])$. Thus $ev$ is smooth along $\pi^{-1}([s])$ by \cite[III. Proposition 10.4.]{Har77}.
\qed}

{\prop \label{qcoker} Let $M$ be a type 2 component and $[s]$ be a general member of $M$. Then the cokernel of the following composition:
\[
T_{M_{red},[s]} \to T_{\overline{M}_{1,0}(X), [s]} = H^0(E, N_s) \to H^0(E, (N_s)_{tf})
\]
has dimension at most $1$. 
}

{\proof  Consider the short exact sequence
\[
0\to T_E \to s^*T_X \to N_s \to 0.
\]
Since $M$ is a type 2 component, we have $h^1(E, s^*T_X) = 1$ and $h^0(E, N_s) = 0$. In particular, the map $h^0(E,s^*T_X) \to h^0(E, N_s)$ is a surjection and the codimension of $T_{M_{red},[s]}$ in $T_{M,[s]}$ is at most $1$. Hence the conclusion follows.
\qed}

\smallskip
We will need to bound the slope of the torsion free part of the normal sheaf. We begin by the following lemma.

{\lemma \label{Lem:genus1bal Sec4} Let $M$ be a type 2 component. Suppose a general member $[f]$ of $M$ satisfies 
\[
\mu^{\max}((N_f)_{tf}) \geqslant \mu^{\min}((N_f)_{tf}) + 3.
\]
Let $a = \mu^{\min}((N_f)_{tf})$ and $p_1,\dots, p_a\in X$ be general points lying on the image of $f$. Then the variety $S$ swept out by deformations of $f$ going through $p_1,\dots, p_a$ is independent of the choices of the $p_i$'s. }

{\proof We may first replace $M$ by its open sublocus parametrizing maps $[g]$ with $\mu^{\min}((N_g)_{tf}) = a$, then replace it again by its image $W$ inside the Kontsevich moduli space $\overline{M}_{1,0}(X)$ and work with the latter. 

Let $W(p_1,\dots,p_a) \subset W$ denote the locus containing $[f]$ parametrizing curves going through the $p_i$'s and let $\pi: U(p_1,\dots, p_a) \to W(p_1,\dots, p_a)$ be its universal family. Let $S = S(p_1,\dots, p_a)$ denote the closure of the variety swept out by $W(p_1,\dots,p_a)$ in $X$ and let $S'$ denote its normalization. By Proposition \ref{twisting the strata of kontsevich space does not dominate}, $S$ has dimension at most $2$, and it suffices to consider the case when $S$ is a surface. Since $U(p_1,\dots, p_a)$ is smooth, the evaluation map $v: U(p_1,\dots, p_a) \to S$ lifts to a morphism $v': U(p_1,\dots, p_a)\to S'$.

Let $q\in S'$ be a point different from the $p_i$'s and consider $Z = v'^{-1}(q)$. By \cite[Lemma 2.4]{She12}, the surface $S'$ is smooth at $q$ and the morphism $v'$ is smooth above $q$, so $Z$ is smooth. Replace $Z$ by its component containing $q$ above $[f]$. We have $Z \cong \pi(Z)$. Hence $Z \cong W(p_1,\dots,p_a, q)$. Let $U'$ denote the universal family over $\pi(Z)$ and $w:U'\to X$ the evaluation map. Since $Z$ is smooth, $U'$ is also smooth, and because $w(U')\subset S$, $w$ lifts to $w': U'\to S'$ as well. We have the dimension count
\begin{align*}
\dim Z &= \dim U(p_1,\dots, p_a) - \dim S'\\
        & = \dim W(p_1, \dots, p_a) + 1 - 2\\
        &\geqslant H^0(N_f(-\sum_{i=1}^a p_i)) - 1 + 1 - 2 \geqslant 1,  
\end{align*}
where the last inequality holds because $H^0(N_f(-\sum_{i=1}^a p_i)) \geqslant 3$. This implies that the family parametrized by  $\pi(Z)$ sweeps out a surface $S(p_1,\dots,p_a,q)$ contained in $S$, and therefore $S(p_1,\dots,p_a,q) = S$. Hence if we start with the set of points $\{ p_1,\dots, p_{a-1}, q\}$, then adding $p_a$ to this set would also give the same surface, i.e.
\[
S(p_1,\dots,p_{a-1}, p_a) = S(p_1,\dots,p_{a-1} ,p_a, q) = S(p_1,\dots,p_{a-1}, q).
\]
By induction, we can replace $p_1,\dots,p_a$ with another set of general points such that the resulting surface $S$ is the same.
\qed}

{\prop\label{mu_min at least 4} Let $M$ be a type 2 component parametrizing maps of degree $e \geqslant 6$. Then for a general member $[f]$ of $M$, we have $\mu^{\min}((N_f)_{tf}) > 3$.}

{\proof
Write $(N_f)_{tf} = N_1\oplus N_2$ where $\deg N_1 = a$ and $\deg N_2 = b$. Suppose to the contrary that $a < 4$. Then since $e\geqslant 6$, Proposition \ref{N_s-tor} implies that $e < b - 1$. Let $p_1,\dots,p_a \in X$ be general points and let $[f]\in M$ be a map passing through the $p_i$'s. Denote by $M' \subset M$ the irreducible component containing $[f]$ and parametrizing maps through the $p_i$'s. 

Since $N_f(-\sum_{i=1}^a p_i)$ is not generically globally generated, the universal family of $M'$ does not dominate $X$ by Proposition \ref{twisting the strata of kontsevich space does not dominate}. We may assume the family of curves parametrized by $M'$ sweeps out a surface $S$. By Lemma \ref{Lem:genus1bal Sec4}, $S$ is independent of the set of points chosen. We have 
\[
    \dim M' \geqslant a + h^0(E, N_f(-\sum_{i=1}^a p_i)) - 1 + 1 = a + b - a - 1 + 1 = b.
\]
In particular, the quantity $h^0(E, N_f(-\sum_{i=1}^a p_i)) - 1$ is the dimension of the family of curves through $a$ general points and the plus one following it accounts for the one-dimensional automorphism of $E$. Then Proposition \ref{bounding the a-inv of a surface} gives
\[
a(S, H|_S) \geqslant \frac{b-1}{e}.
\]
Since $e < b-1$, we have 
\[
a(S, H|_S) > \frac{b-1}{b-1} = 1.
\]
Hence $a(S,L) = 2$ by Theorem \ref{Classifya-inv} and $M$ parametrizes covers of lines by Corollary \ref{facunifamlinecor}. A contradiction. Thus $\mu^{\min}((N_f)_{tf}) = a > 3$.

{\theorem \label{no type 2}  There is no type 2 component parametrizing maps of degree at least $6$.
}

{\proof Suppose to the contrary that there is such a type 2 component $M$. Let $W$ be a variety equipped with a generically finite morphism to the locus induced by $M$ in $\overline{M}_{1,0}(X)$. By Theorem \ref{nfconn} and Proposition \ref{evflat}, a general fibre of the evaluation map $ev_W: \CU_W \to X$ is irreducible and contained in the flat locus. Hence we are in the situation of Grauert-M\"ulich (Theorem \ref{GM}). 

Write the Harder Narasimhan filtration of $s^*T_X$ as 
\[
0 = \CF_0 \subset \CF_1 \subset \dots \subset \CF_k = s^*T_X.
\]
Denote by $t$ the length of the torsion part of $N_s$, by $\CG$ the subsheaf of $(N_s)_{tf}$ generated by global sections, by $V$ the tangent space to $W$ at $s$, and by $q$ the dimension of the cokernel of the composition
\[
V \to T_{\overline{M}_{1,0}(X),s} = H^0(E,N_s) \to H^0(E, (N_s)_{tf}).
\]
Theorem \ref{GM} shows that for every index $1\leqslant i \leqslant k-1$, we have
\[
\mu(\CF_i/\CF_{i-1}) - \mu(\CF_{i+1}/\CF_i) \leqslant (q+1)\mu^{max}({M_\CG\smvee}) + t.
\]
By assumption, $(N_s)_{tf}$ is globally generated, so by \cite[1.3 Corollary.]{But94} and Proposition \ref{mu_min at least 4}, we have 
\[
\mu^{max}({M_\CG\smvee}) \leqslant \frac{\mu^{\min}(\CG)}{\mu^{\min}(\CG) - 1} = 1 + \frac{1}{\mu^{\min}(\CG) - 1} \leqslant 1 + \frac{1}{3} < \frac{3}{2}.
\]
Since $t \leqslant 1$ by Proposition \ref{N_s-tor} and $q \leqslant 1$ by Proposition \ref{qcoker}, we have
\[
\mu(\CF_i/\CF_{i-1}) - \mu(\CF_{i+1}/\CF_i) < (1+1)\cdot\frac{3}{2} + 1 = 4.
\]
By assumption, $s^*T_X$ contains $\CO_E$ as a factor. The slope panel $SP(s^*T_X) = (s_1, s_2, s_3) = (s_1, s_2, 0)$. Since $s$ has degree at least $6$, $s_1+s_2 \geqslant 12$. But this contradicts the bound on the difference between the slope of successive quotients in the filtration. Hence there is no type $2$ component parametrizing maps of degree at least $6$.
\qed}

\smallskip
Finally, a similar argument to Theorem \ref{no type 1} rules out the type $3$ components.

{\theorem \label{no type 3} There is no type 3 component parametrizing maps of degree at least $4$. } 

{\proof Suppose to the contrary that there is such a type 3 component $M$. Let $[s]$ be a general member of $M$ with $s^*T_X = \CO_E\oplus \CO_E\oplus \CL$. Let $e$ denote the degree of $s$. Let $p\in E$ be a general point and let $W$ be an irreducible component of the locus of curves containing $[s]$ through $s(p)$. Since $s^*T_X(-p)$ is not generically globally generated, the universal family of $W$ sweeps out a surface $S$ in $X$ by Proposition \ref{nondom}. We have
\[
\dim W \geqslant h^0(E, s^*T_X(-p)) - h^1(E, s^*T_X(-p)) + 1 = 2e - 3 + 1 = 2e - 2,
\]
where the plus one accounts for varying the choice of $p$ on $E$. Hence
\[
a(S,L) \geqslant \frac{2e - 3}{e}
\]
by Proposition \ref{bounding the a-inv of a surface}. Then $a(S,L) > 1$ if $e\geqslant 4$. But this contradicts Theorem \ref{Classifya-inv} and Corollary \ref{facunifamlinecor}.
\qed}

\smallskip
Gathering all the pieces together, we obtain the following theorem.

{\theorem \label{no non-free-comp-morphism} There is no non-free component in $\Mor(E,X)$ generically parametrizing a dominant family of maps of degree $e\geqslant 6$ that are birational onto their images. 
}

\smallskip
This allows us to deduce that any irreducible main component of the Kontsevich moduli space parametrizes stable maps whose domain curves vary in moduli.

{\corollary \label{Cor: mor does not dominate kontsevich} Let $M$ be an irreducible component of $\Mor(E,X)$. Then under the natural morphism $\phi: \Mor(E,X) \to \overline{M}_{1,0}(X)$, the image of $M$ does not dominate any irreducible component of $\overline{M}_{1,0}(X)$.
}
 
{\proof By Theorem \ref{irred of N_e}, Theorem \ref{irred of double covers of conics}, and Theorem \ref{no non-free non-dominant comp}, we may assume $M$ parametrizes a dominant family of morphisms that are birational onto their images. Let $e$ be the degree of the morphisms $M$ parametrizes. 

If $e\leqslant 5$, then Corollary \ref{base1} implies that there is a unique irreducible component $\overline{M}_e$ in $\widetilde{M}_{1,0}(X,e)$ parametrizing stable maps that are birational onto their images. Denote by $\overline{M}'_e$ the universal family of $\overline{M}_e$. Since $\overline{M}_e$ is birational to $\overline{\CH}^{e,1}$ and $\overline{\CH}^{e,1}$ parametrizes a reducible genus one curve whose stabilization is the nodal rational curve $[C] \in \overline{M}_{1,0}$, the map $\overline{M}'_e \to \overline{M}_{1,1}$ is surjective. Hence $M$ does not dominate $\overline{M}_e$ under $\phi$.

If $e\geqslant 6$, then by Theorem \ref{no non-free-comp-morphism}, we may assume a general member $[f]$ of $M$ has $h^1(E, f^*T_X) = 0$. Since $T_{M, [f]} = H^0(E, f^*T_X)$ and the map $H^0(E, f^*T_X) \to H^0(E, N_f)$ is not a surjection, the statement follows.
\qed}

\section{\label{Free components} Free components}

In the last section, we first proceed to prove the Movable Bend-and-Break Theorem \ref{Thm:MBB}. We then use this result together with the irreducibility of low degree components to inductively prove the irreducibility of free birational component in each degree.

Recall that in our Definition \ref{free}, a morphism $f:E\to X$ is a free curve if $f^*T_X$ is generically globally generated and $h^1(E, f^*T_X) = 0$. A free curve is a smooth point in both $\Mor(E,X)$ and $\overline{M}_{1,0}(X)$. However, the vanishing of $h^1(E, f^*T_X)$ is a stronger condition than the vanishing of $h^1(E, N_f)$. Since we will mostly be working with the Kontsevich space, in the rest of the section, we turn our attention to the normal sheaf. In particular, we apply Grauert-M\"ulich to deduce when the normal sheaf is locally free for a general map.

We first show that for a free component, a general fibre of the evaluation map is irreducible.

{\prop Let $M \subset \Mor(E,X)$ be a free component parametrizing maps of degree $e\geqslant 4$. Let $U$ be the universal family of $M$. Then a general fibre of the evaluation map $ev: U \to X$ is irreducible.}

{\proof Suppose to the contrary that a general fibre of the evaluation map is reducible. We take the Stein factorization of $ev:\widetilde{U} \to X$, where $\widetilde{U}$ is a resolution of singularities. We let $Y$ denote a resolution of the finite part of the Stein factorization. This produces a generically finite morphism $f: Y \to X$, and $M$ gives rise to a dominant component $N$ of $\Mor(E, Y)$. Denote the induced morphism as $\phi: N \to M$. We have $\dim M = \dim N$.

Let $s': E\to Y$ be a general member of $N$ such that $s = f\circ s'$ is a general member of $M$. Let $y\in Y$ be a general point in $Y$. Denote the family of curves parametrized by $N$ through $y$ by $N_y$ and the family of curves parametrized by $M$ through $x = f(y)$ by $M_x$. We claim that $M_x$ parametrizes a dominant family of curves on $X$. Suppose to the contrary and let $S$ be the surface swept out by the family of $M_x$. We have
\[
\dim M_X \geqslant h^0(E, s^*T_X(-p)) - h^1(E, s^*T_X(-p)) + 1 = 2e - 3 + 1 = 2e - 2,
\]
where the plus one encodes the choice of varying the point $p\in E$ sending to $x$. By Proposition \ref{bounding the a-inv of a surface}, we have
\[
a(S, H|_S) \geqslant \frac{\dim M_S - 1}{e} \geqslant 2 - \frac{3}{e}.
\]
This implies that if $e \geqslant 4$, then $a(S, H|_S) > 1$, and $M_x$ parametrizes covers of lines by Theorem \ref{Classifya-inv} and Corollary \ref{facunifamlinecor}. A contradiction.

This implies that $N_y$ also parametrizes a dominant family of curves on $Y$. In particular, we have $h^1(E, s'^*T_Y) \leqslant 1$ because otherwise $s'^*T_Y$ must contain $\CO_E\oplus\CO_E$ as a direct summand and the universal family of $N_y$ does not dominate $Y$ by Proposition \ref{nondom}. 

Consider the following equality
\[
-K_X\cdot s_*E = \dim M = \dim N \leqslant -K_Y\cdot s'_*E + h^1(E, s'^*T_Y) = -K_Y\cdot s'_*E + 1.
\]
We can rewrite the inequality in the following way
\[
\left(K_Y + (\frac{3}{2} + \epsilon) f^*H + (\frac{1}{2} - \epsilon)f^*H\right)\cdot s'_*E \leqslant 1,
\]
where $\epsilon > 0$ is a sufficiently small constant. We have
\[
\left(K_Y + (\frac{3}{2} + \epsilon) f^*H\right)\cdot s'_*E \leqslant 1 - (\frac{1}{2} - \epsilon)f^*H\cdot s'_*E.
\]
Since $e \geqslant 4$, the right hand side of the inequality is less than $0$. Hence $a(Y, f^*H) \geqslant \frac{3}{2} + \epsilon$, and Theorem \ref{a-inv of 3fold} implies that $a(Y, f^*H) \geqslant 2$. This contradicts Theorem \ref{Classifya-inv} and Corollary \ref{facunifamlinecor}.
\qed}

{\prop \label{locfreegenusone} Let $M \subset \Mor(E,X)$ be a free component parametrizing maps that are birational onto their images. A general member of $M$ has locally free normal sheaf. }

{\proof  If $M$ parametrizes maps of degree at most $5$, then the assertion follows from Theorem \ref{basecase}. Hence we may assume that $M$ parametrizes maps of degree at least $6$. Let $f$ be a general member of $M$. \cite[1.8 Theorem]{Kol96} implies that if $f$ is not an immersion, then $h^1(E, f^*T_X(-2p)) > 1$ for some $p\in E$. Hence $f^*T_X$ must be decomposable and we may write $f^*T_X = \CF \oplus \CG$ with $h^1(E, \CF(-2p)) > 1$ and $h^1(E, \CG(-2p)) = 0$. Note that $\CF \neq f^*T_X$ because $\deg \CF \leqslant 2\rank(\CF) \leqslant 6 < \deg f^*T_X$. 

If $\rank(\CF) = 1$, then $h^1(E, \CF(-2p)) > 1$ implies that $\deg \CF \leqslant 0$. This contradicts the fact that $f$ is a free curve. Hence $\rank(\CF) = 2$ and $\CF$ is either indecomposable of degree at most $2$, or decomposable of the following forms: 
\[
\CO_E(2p) \oplus \CO_E(2p),\  \CO_E(q) \oplus \CO_E(2p),\  \CO_E(q) \oplus \CO_E(p),
\] 
where $p$ and $q$ are possibly the same point on $E$. Let the slope panel of $f^*T_X$ be $SP(f^*T_X) = (a_1, a_2, a_3)$. It is easy to check that we are in the situation of Grauert-M\"ulich \ref{GM}. In particular, we have $a_i - a_{i+1} \leqslant 5$ for $i = 1, 2$. Hence $a_1 + a_2 + a_3 < 12$ for the above possible choices of $\CF$. This contradicts the fact that $\deg f^*T_X \geqslant 12$. 
\qed}

\smallskip
We also include an analogous statement for the space of rational curves. The proof is omitted as it is similar to the previous proposition.

{\prop \label{veryfreerational} For $e\geqslant 2$, the irreducible component $M\subset \Mor(\BP^1, X)$ parametrizing maps of degree $e$ that are birational onto their images is very free.}

\smallskip
Next, we show that the normal bundle will be balanced for a general member of a free component.

{\definition Let $f:E\to X$ be an immersion of degree $e$. We say $N_f$ is balanced if the slope panel satisfies $SP(N_f) = (e, e)$. Otherwise $N_f$ is unbalanced. We say an irreducible component $M\subset \widetilde{M}_{1,0}(X)$ is balanced (resp. unbalanced) if $N_f$ is balanced (resp. unbalanced) for a general $[f]\in M$.}

{\prop \label{bal}  Let $M$ be a free component in $\widetilde{M}^{bir}_{1,0}(X,e)$. Then $M$ is a balanced component.}

{\proof By Proposition \ref{locfreegenusone}, a general member of $M$ has locally free normal sheaf. Let $[f]$ be a general member of $M$. Suppose $N_f$ is unbalanced. Then we can write $N_f = N_1\oplus N_2$ where $\deg N_1 = a$, $\deg N_2 = b$, and $b-a\geqslant 2$. Let $p_1,\dots,p_a \in X$ be general and let $[f]\in M$ be a map passing through these general points. Denote by $M' \subset M$ the irreducible component containing $[f]$ and parametrizing maps passing through $p_1,\dots, p_a$. 

Since $N_f(-\sum_{i=1}^a p_i)$ is not generically globally generated, $M'$ parametrizes a family of curves on a surface $S\subset X$. Similar argument to Lemma \ref{Lem:genus1bal Sec4} shows that $S$ is independent of the set of points chosen. Replace $S$ by a resolution of singularities. Then $M'$ gives rise to a component $M" \subset \overline{M}_{1,0}(S)$. Since there is a dominant family of curves going through $a$ general points on $S$ and this family has dimension $h^0(E, N_f(-\sum_{i=1}^a p_i))$, we can compute the lower bound of the dimension of $M"$: 
\[
\dim M" \geqslant a + h^0(E, N_f(-\sum_{i=1}^a p_i)) = a + b - a = b.
\]
Write $C$ for the induced image of $f$ on $S$, and denote by $L$ the pullback of the hyperplane section of $X$ to $S$. We have
{\allowdisplaybreaks\begin{align*}
(K_{S} + a(S, L)L)\cdot C &= K_{S} \cdot C + a(S,L)L\cdot C\\
					&= - \dim M" + a(S,L)\cdot e\\
					&\leqslant -b + a(S,L)\cdot e\\
					&< -e + a(S,L)\cdot e = e(a(S,L) - 1).
\end{align*}}

Since $K_{S} + a(S, L)L$ is pseudo-effective, we have $a(S,L) > 1$. Hence by Theorem \ref{Classifya-inv}, $a(S,L) = 2$ and $S$ factors through the universal family of curves. This implies that $M"$, and hence $M$, generically parametrizes covers of lines, a contradiction to Corollary \ref{facunifamlinecor}. Hence $N_f$ is balanced. 
\qed
}

\subsection{Towards Movable Bend-and-Break} \hfill

\smallskip
\noindent From now on, we will mainly be interested in the components of  $\widetilde{M}_{g,0}(X)$ that generically parametrize maps that are birational onto their images. We denote the union of such components by $\widetilde{M}^{bir}_{g,0}(X)$. 

In this subsection, we prove the Movable Bend-and-Break Lemma \ref{Prop:MBB} by studying the deformation of stable maps through fixed general points and general curves. We first explain what we mean by ``general curves":

{\definition \cite[Section 6.2]{BLRT22} Let $p:U\to R$ be a family of irreducible and reduced curves with evaluation map $ev: U\to X$. We say the family $p$ is basepoint free if the evaluation map is flat.
}

\smallskip
Basepoint free family of curves can be obtained by intersecting general members of a very ample linear series. A key property of such families is that given any fixed codimension two subset in $X$, a general member will not intersect this subset. The following construction will be used many times in the rest of this section:

{\construction \label{construction of T} Let $M\subset \widetilde{M}^{bir}_{1,0}(X,e)$ be a free component. Let $M^\circ$ be the open subset of $M$ parametrizing free curves. There are two possible scenarios for a general member in $M^\circ$:
\begin{enumerate}
\item The normal bundle is indecomposable of degree $2e$ or is a direct sum of line bundles of the same degree.

\item The normal bundle is of the form $\CL_1\oplus \CL_2$ where $l_1 = \deg \CL_1$, $l_2 = \deg \CL_2$ and $0 < l_1 < l_2$.
\end{enumerate}
Then for any positive integers $r$ and $s$ satisfying (in the respective cases)
\begin{enumerate}
\item $r \leqslant e - 1$ and $s \leqslant 2e - 2r - 1$, 

\item $r \leqslant l_1$ and $s \leqslant l_2 - l_1$,
\end{enumerate}
denote by $T_{r,s}$ the locus in $M^\circ$ parametrizing curves going through $r$ general points and $s$ general members of $p$.}

{\lemma \cite[c.f. Lemma 6.4]{BLRT22} \label{IncidenceCorr}  Let $X$ be a smooth projective threefold. Let $p: U \to R$ be a basepoint free family of curves on $X$. Let $M\subset \widetilde{M}^{bir}_{1,0}(X,e)$ be a free component and let $T_{r,s}$ be the locus constructed in Construction \ref{construction of T}. Then $T_{r,s}$ satisfies
\[
\codim(T_{r,s}) \geqslant 2r+s.
\]
If furthermore $p$ is constructed by taking complete intersections of a big and basepoint free linear series on $X$, then the above achieves equality.
}

\smallskip
Next, we address the local behavior of a stable map with a contracted genus one component. Such a stable map may be a singular point in the moduli space and could lie in the intersection of multiple components in $\overline{M}_{1,0}(X)$.

{\lemma \cite[5.9. Lemma.]{Vak00} \label{smoothable} Let $Z$ be a complete reduced nodal curve of arithmetic genus 1, and let $\pi: Z\to X$ be a morphism to a smooth cubic threefold. Assume $(Z, \pi)$ can be smoothed. If $B$ is a connected union of contracted components of $Z$ of arithmetic genus 1, intersecting $\overline{Z\backslash B}$ in $k$ points, and $t_1, \dots, t_k$ are tangent vectors to $\overline{Z\backslash B}$ at those points, then $\{\pi(t_i)\}_{i=1}^k$ are linearly dependent in $T_{X, {\pi(B)}}$.
}

\smallskip
When $k = 1$ in the above lemma, the domain of the stable map consists of at least two irreducible components, of which the genus one curve is contracted. In particular, the image has a cuspidal singularity at the contracted image point. We compute the codimension of the locus of curves with cuspidal singularity in $\overline{M}_{0,0}(X)$:

{\lemma \label{cusp} Let $M \subset \overline{M}_{0,0}(X)$ be a free component generically parametrizing maps birational onto their images of degree $e\geqslant 4$. The locus of $M$ parametrizing very free curves whose images have cuspidal singularities has codimension at least $2$ in $M$. 
}

{\proof By Proposition \ref{veryfreerational}, the assumption on degree implies that $M$ is a very free component. Since the map $\Mor(\BP^1, X) \to \overline{M}_{0,0}(X)$ is flat over $M$, it suffices to show the equivalent statement for the morphism space. 

Denote by $M'$ the corresponding component in $\Mor(\BP^1, X)$. Let $U\subset M'$ be the locus of very free curves. Since the locus of maps whose images have cuspidal singularities is a closed subscheme, it suffices to look at the deformations of a general map $[f] \in U$ with a cuspidal singularity at $p\in \BP^1$. We estimate the dimension of the following locus:
\[
W \coloneqq \{[g] \in U \ |\ g \ \text{is not an immersion}\},
\] 
where $[g]$ denotes a deformation of $[f]$. Let $x = g(p)$. Let $J \cong \Spec k[\epsilon]/\epsilon^2$ denote the length $2$ subscheme representing the point $p$ and its tangent direction $t_p$. Let $g_{p} = g|_J$. By varying $p\in \BP^1$ and $x\in X$, we can write $W$ set-theoretically as the following:
\[
W = \bigcup_{p\in\BP^1} \bigcup_{x\in X} \Mor(\BP^1, X\ |\ g_{p}) \cap U.
\]
By assumption, $U$ parametrizes very free curves, so we have

\begin{align*}
    \dim \Mor(\BP^1, X\ |\ g_{p}) \cap U &=  h^0(\BP^1, f^*T_X(-2p)) \\
        &= h^0(\BP^1, f^*T_X) - 6\\
        &= \dim M' - 6.
\end{align*}

Hence, we have $\dim W \leqslant \dim M' - 2$, i.e. the codimension of the locus of very free curves with cuspidal singularity is at least $2$ in $M'$.
\qed}

\smallskip
The following lemma shows that a general free curve intersects a codimension two subscheme of $X$ transversally. 

{\lemma \label{tangent intersection has high codimension}
Let $C \subset X$ be a curve. Let $M \subset \overline{M}_{0,0}(X,e)$ be a very free component. Then the locus of degree $e$ very free curves in $M$ intersecting $C$ tangentially has codimension at least $3$.}

{\proof Since the map $\Mor(\BP^1, X) \to \overline{M}_{0,0}(X)$ is flat over $M$, it suffices to show the equivalent statement for the morphism space. Given a point $p\in\BP^1$, denote by $J \cong \Spec k[\epsilon]/\epsilon^2$ the length $2$ subscheme representing the point $p$ and its tangent direction $t_p$. Let $g_{v,p}: J\to X$ denote the map sending $p$ to $x$ and $t_p$ to $v\in T_{X,x}$. Let $U\subset \Mor(\BP^1,X,e)$ denote the irreducible locus of very free rational curves. Let $W\subset U$ be the sublocus of curves intersecting $C$ tangentially. We can write $W$ set-theoretically as the following:
\[
W = \bigcup_{x\in C} \bigcup_{v\in T_{C,x}}\bigcup_{p\in \BP^1} \Mor(\BP^1, X\ |\ g_{p,v}) \cap U.
\]
Since $U$ is the locus of very free rational curves, we have
\[
\dim \Mor(\BP^1, X\ |\ g_{p,v}) \cap U = h^0(\BP^1, f^*T_X(-2p)) = h^0(\BP^1, f^*T_X) - 6.
\]
Thus the dimension of $W$ is at most
\[
\dim W \leqslant h^0(\BP^1, f^*T_X) - 6 + 3 = \dim U - 3.
\]
\qed}

\smallskip We are ready to prove the Movable Bend-and-Break Lemma. Proposition \ref{bal} implies that every component $M$ in $\widetilde{M}_{1,0}^{bir}(X,e)$ generically parametrizing free curves is a balanced component. Hence given a general map $[f]\in M$ with balanced normal bundle, we may construct a $1$-dimensional locus in $M$ parametrizing deformations of $f$ through $e - 1$ general points and one general member of a basepoint free family of curves by Lemma \ref{IncidenceCorr}.

{\prop[Movable Bend-and-Break Lemma] \label{Prop:MBB}  Let $p:U\to R$ be a basepoint free family of curves on $X$. Let $M$ be an irreducible component in $\widetilde{M}_{1,0}^{bir}(X,e)$ with $e\geqslant 4$ and let $M^\circ$ be the open subset parametrizing free curves with irreducible domains. Let $T \subset M$ be the closure of the 1-dimensional locus $T_{e-1,1}$ as constructed in Construction \ref{construction of T}. Suppose $T$ parametrizes a stable map $f:Z\to X$ such that $Z$ contains a smooth genus one component and a rational component. Then $Z$ consists of exactly two irreducible components $Z_0$ and $Z_1$ satisfying the following conditions: 
\begin{enumerate}
    \item $Z_0$ is a free smooth genus one curve.
    \item $Z_1$ is a free rational curve.
    \item $f|_{Z_i}$'s are birational onto their images.
\end{enumerate} 
}

{\remark \label{Remark: MBB} We explain what we mean by ``general" in the above proposition. The locus $W$ swept out by non-free rational curves and non-dominant genus one curves of degree at most $e$ in $X$ is a proper subset. We remove this locus from our choice of general points. Then we choose the points and the curve of the basepoint free family to be general so that they impose independent conditions on the families of free rational curves and free genus one curves of degree at most $e$. In particular, for degree $3\leqslant d < e$, we may ensure that there is a finite set of curves of degree $d$ through $d$ general points such that every member is free. We may also guarantee that there is a $1$-dimensional family of curves of degree $d$ through $d-1$ general points and a general basepoint free curve, and that a general member of this family is free. The proof of both are similar and we give one for the former in the genus one case.

Let $T = T_{d-1,0} \subset \widetilde{M}^{bir}_{1,0}(X,d)$ denote the $2$-dimensional locus parametrizing stable maps through general points $p_1,\dots, p_{d-1}$ built in Construction \ref{construction of T}. Denote by $\CU_T$ the universal family of $T$ and by $ev: \CU_T \to X$ the evaluation map. Since the locus of non-free genus one curves in $\widetilde{M}^{bir}_{1,0}(X,d)$ has codimension at least one by Corollary \ref{Cor: mor does not dominate kontsevich}, the family of non-free curves parametrized by $T$ has dimension at most one. Hence the images of the non-free curves parametrized by $T$ under $ev$ has codimension at least $1$ in $X$. Similarly, the locus swept out by reducible curves in $T$ also has codimension at least one in $X$. Hence we may pick a general point $p_d$ so that the locus $T_{d,0}$ constructed from $p_1,\dots,p_d$ is a finite set of free genus one curves.
}

{\proof 
Let $p_1,\dots,p_{e-1}$ be the general points and let $Q$ be the general member of $p$ used in the construction of $T_{e-1,1}$. Write
\[
Z = Z_0 + \sum_{i = 1}^r Z_i + \sum_{i=r+1}^k Z_i + Z_c,
\]
where $Z_0$ is a smooth genus one curve, $Z_i$ are free rational curves for $1\leqslant i \leqslant r$, non-free rational curves for $r+1\leqslant i \leqslant k$, and $Z_c$ is the sum of contracted rational components. Let $e_i$ denote the $H$-degree of $Z_i$ and let $a_i$ denote the number of general points $Z_i$ goes through for $i = 0,\dots, r$. Denote by $e'$ the total $H$-degree of non-free rational curves. We divide our discussion based on the image of $Z_0$.

Suppose that $Z_0$ is contracted to a point. In this case, $a_0 \leqslant 1$. Assume first that $a_0 = 0$ so that $f(Z_0)$ is not one of the general points. Since $Z$ goes through $e-1$ general points, we have the following inequality:
\[
e - 1 \leqslant \sum_{i=1}^r a_i \leqslant \sum_{i=1}^r e_i = e - e'.
\]
Hence $e' \leqslant 1$. If $e' = 1$, then there is a non-free line $L$ and we have $a_i = e_i$ for all $i$, so each free rational curve goes through the maximum number of points. These points determine a finite set of possibilities for each $f|_{Z_i}$, so they are all disjoint. Since the image of $Z$ is connected, the non-free line $L$ must meet the $Z_i$'s and the basepoint free curve $Q$. On the other hand, the family of non-free lines is $1$-dimensional, and by a dimension count, $L$ can only meet at most one of the $Z_i$'s and $Q$. Hence $f(Z)$ cannot be connected.

If $e' = 0$, all but one of the free rational curves must contain $a_i$ general points and have degree $e_i = a_i$. The remaining component, say $Z_1$, must contain $a_1$ general points and have degree $e_1 = a_1 + 1$. In particular, the general points determine a finite set of choices of $f|_{Z_i}$'s for each $i>1$. Since $Q$ is general, it does not intersect a codimension $2$ locus in $X$, so $f(Z_1)$ must meet $Q$. Hence $Z_1$ deforms in a family of dimension at most $1$. Since it is a codimension one condition for a family of free curves to meet a fixed curve, the generality of the points $p_i$'s  ensures that there can be at most two free rational components for the image of $Z$ to be connected. By Lemma \ref{smoothable}, either $Z_1$ is the only rational curve with a cuspidal singularity at $f(Z_0)$, or $Z_1$ and $Z_2$ meet tangentially at $f(Z_0)$. By Lemma \ref{cusp}, the locus of very free cuspidal rational curves in this degree range has codimension at least $2$, so the former case is not possible. The latter case is also not possible because Lemma \ref{tangent intersection has high codimension} shows that $Z_1$ and $Z_2$ intersect transversally.

Assume now that $a_0 = 1$ so that $f(Z_0)$ is one of the general points. Let $b$ be the number of free rational curves through $f(Z_0)$. Then we have the inequality
\[
(e - 1) + b - 1 \leqslant \sum_{i=1}^r a_i \leqslant  \sum_{i=1}^r e_i = e - e',
\]
where the $b - 1$ on the left hand side indicates that $f(Z_0)$ is counted an extra $b - 1$ times. This shows that $b = 1$ or $b = 2$. 

If $b = 2$, then there are two free rational curves, say $Z_1$ and $Z_2$ going through $f(Z_0)$ and the inequality becomes
\[
e \leqslant \sum_{i=1}^r a_i \leqslant  \sum_{i=1}^r e_i = e - e'.
\]
This implies that there is no non-free rational curves and the free rational curves go through the maximum number of general points and hence only have finitely many choices. But this implies that they are disjoint from $Q$.

Hence $b = 1$ and the above inequality becomes
\[
e - 1 \leqslant \sum_{i=1}^r a_i \leqslant  \sum_{i=1}^r e_i = e - e'. 
\]
If $e' = 1$, there is a non-free line $L$. In this case, $L$ must intersect $Q$, and the image of $Z$ cannot be connected. If $e' = 0$, then there must be exactly one free rational curve, say $Z_1$, that goes through $a_1$ general points and has degree $e_1 = a_1 + 1$. In particular, $Z_1$ must meet $Q$. Hence $Z_1$ deforms in a family of dimension at most $1$. Again by Lemma \ref{smoothable}, $Z_1$ has a cuspidal singularity at $f(Z_0)$, which is not possible because Lemma \ref{cusp} shows that the locus of very free cuspidal rational curves has codimension at least $2$. 

Hence, $Z_0$ cannot be contracted. We have the inequality:
\[
e - 1 \leqslant \sum_{i=0}^r a_i \leqslant \sum_{i=0}^r e_i = e - e'.
\]
This implies $e' = 1$ or $e' = 0$. The first case is not possible by similar arguments. Thus, there is no non-free rational curve. Moreover, there exists a single $0\leqslant j \leqslant r$ such that $Z_j$ goes through $e_j - 1$ general points and all other $f|_{Z_i}$'s for $i\neq j$ have finitely many possibilities upon meeting the maximum number of general points. Hence $f|_{Z_j}$ must meet $Q$ and it deforms in a $1$-dimensional family. Since the image of $Z$ is connected, this implies that there are exactly two components for similar reasons as previous.  

We may write $Z = Z_0 + Z_1$. The genus one component $Z_0$ has degree $e_0$ and goes through at least $e_0-1$ general points. If $f|_{Z_0}$ is a cover of some smaller degree curve $C$, then $\deg C \leqslant \frac{e_0}{2}$. Since $C$ goes through at least $e_0-1$ general points, we must have 
\[
2\cdot (e_0-1) \leqslant e_0,
\]
which implies that $e_0\leqslant 2$, so $f|_{Z_0}: Z_0\to X$ is a degree 2 cover of a free line. Since there are only finitely many possible lines through a general point, $f|_{Z_0}$ does not meet $Q$. On the other hand, all other rational components go through the maximum number of general points, and hence they cannot meet $Q$. A contradiction. Consequently,  $f|_{Z_0}$ is birational onto its image. The same argument shows that $f|_{Z_1}$ is also birational onto its image. Finally, we show that $f|_{Z_0}$ is free. If $Z_0$ goes through $e_0$ general points, then the claim follows from Remark \ref{Remark: MBB}. If $Z_0$ goes through $e_0 - 1$ general points and $Q$, then it deforms in a $1$-dimensional family. Since there is a finite number of non-free curves in this family by Remark \ref{Remark: MBB} again, the choice of general points guarantees that $f|_{Z_1}$ intersects the free curves of this family, and thus $f|_{Z_0}$ is free.
\qed}

\smallskip
The significance of the Movable Bend-and-Break Lemma $\ref{Prop:MBB}$ is that by varying the number of constraints we impose on the deformation, it tells us precisely the codimension of the locus corresponding to different broken types emerging in the locus $T$ constructed.

\subsection{Deforming nodal curves} \hfill

\smallskip
\noindent The limit of a one-parameter family of stable maps with domain a smooth genus one curve may have domain an irreducible nodal rational curve. We investigate how stable maps with nodal rational domains deform in $\widetilde{M}_{1,0}(X)$.

{\definition Given $[f] \in \overline{M}_{0,0}(X)$. We say $f$ is a stable map of nodal rational type if $f$ factors through a rational curve of arithmetic genus one.}

\smallskip
Let $C$ denote a nodal rational curve. Then a stable map $f:C\to X$ in $\widetilde{M}_{1,0}(X)$ can be regarded as a stable map of nodal rational type $\tilde{f}:\widetilde{C}\to X$ by precomposing $f$ with a partial normalization $\widetilde{C} \to C$. This allows us to relate the deformation theory of $f$ with that of $\tilde{f}$.

{\prop \label{free genus one nodal stable map} Let $f: C \to X$ be a genus one stable map such that $C$ is an irreducible nodal rational curve. Denote by $\tilde{f}: \BP^1 \to X$ the composition of $f$ with the normalization $\nu: \BP^1 \to C$. If $\tilde{f}$ is very free, then $f$ is a smooth point of $\widetilde{M}_{1,0}(X)$.}

{\proof By \cite[Section 1.2]{Tes05}, there is an isomorphism between the dual of the obstruction space of $[f]$ and $H^0(C, \CC_f\otimes\omega_C)$, where $\CC_f$ is the conormal sheaf defined in \cite[Definition 1.9]{Tes05} as the cokernel of the exact sequence 
\[
0 \to \CC_f \to f^*\Omega_X \to \Omega_C \to Q_f \to 0.
\]
Since $\omega_C \cong \CO_C$, we have $H^0(C, \CC_f\otimes\omega_C) = H^0(C, \CC_f) = 0$. This gives an injection $H^0(C, \CC_f\otimes\omega_C) \to H^0(C, f^*\Omega_X)$. On the other hand, we have the normalization exact sequence
\[
0 \to \CO_C \to v_*\CO_{\BP^1} \to \CO_p \to 0,
\]
where $p$ is the node of $C$. Tensoring with the locally free sheaf $f^*\Omega_X$ and using the projection formula gives a short exact sequence
\[
0 \to f^*\Omega_X \to \nu_*\tilde{f}^*\Omega_X \to {f^*\Omega_X}|_p \to 0.
\]
Hence there is an injection 
\[
H^0(C, f^*\Omega_X) \to H^0(C, \nu_*\tilde{f}^*\Omega_X) \cong H^0(\BP^1, \tilde{f}^*\Omega_X) \cong H^1(\BP^1, \tilde{f}^*T_X\otimes\CO_{\BP^1}(-2)).
\]
Since $\tilde{f}$ is a very free rational curve, we have $H^1(\BP^1, \tilde{f}^*T_X\otimes\CO_{\BP^1}(-2)) = 0$. This implies that $H^0(C, f^*\Omega_X) = 0$, and hence $H^0(C, \CC_f\otimes\omega_C) = 0$ and $[f]$ is a smooth point.
\qed}

\smallskip
Moreover, we may regard the locus of genus one stable maps with nodal rational domains as a sublocus in $\overline{M}_{0,0}(X)$. The following lemma implies that the locus of stable maps of nodal rational type through general points and general basepoint free members has the expected dimension.

{\lemma \label{dimension of nodal rational curve} 
Let $M\subset \widetilde{M}^{bir}_{1,0}(X,e)$ be a free component with $e\geqslant 3$ and let $T_{r,s}\subset M$ be the locus constructed in Construction \ref{construction of T}. Then the sublocus of stable maps with nodal rational domains in $T_{r,s}$ has dimension equal to $\dim T_{r,s} - 1$. 
}

{\proof We first note that the locus $W$ of stable maps with nodal rational domains in $M$ has codimension one. Indeed, if $W$ is dense, then it induces a dense locus in $\overline{M}^{bir}_{0,0}(X,e)$ since the latter is irreducible by Theorem \ref{Starr}. Moreover, it is a free component by \cite[Corollary 7.7]{LT19a}. Hence Proposition \ref{veryfreerational} and \cite[1.8 Theorem]{Kol96} implies that a general member of $\overline{M}^{bir}_{0,0}(X,e)$ is an embedding. This is a contradiction.

Denote by $T'_{r,s}$ the preimage of $T_{r,s}$ in $\widetilde{M}_{1,1}^{bir}(X,e)$. Then $\dim T'_{r,s} = \dim T_{r,s} + 1$. Since the points and the members of the basepoint family of curves are general, the morphism $\pi: T'_{r,s} \to \overline{M}_{1,1}$ is surjective. Hence the fibre of $\pi$ above the nodal rational curve $[C, p] \in \overline{M}_{1,1}$ has dimension equal to $\dim T_{r,s}$. This shows that the intersection of $T_{r,s}$ with the locus of stable maps with nodal rational domains has dimension equal to $\dim T_{r,s} - 1$. 
\qed}

\smallskip
In particular, for $e\geqslant 3$, the locus $T_{e-1,0} \subset \widetilde{M}_{1,0}^{bir}(X,e)$ contains a $1$-dimensional sublocus $T^C_{e-1,0}$ that parametrizes stable maps with nodal rational domains. By partially normalizing the domains of the stable maps in $T^C_{e-1,0}$ so that they become nodal curves of genus zero, we obtain a morphism: 
\[
T^C_{e-1,0} \to \overline{M}^{bir}_{0,0}(X,e).
\]
Hence the locus of stable maps of nodal rational type through $e-1$ general points has dimension one. We will apply a nodal version of the Movable Bend-and-Break to this locus. Before stating the result, we need the following definition:

{\definition Let $C$ be a genus one nodal curve. We say $C$ is a banana curve if $C$ consists of exactly two irreducible smooth rational curves intersecting each other at two distinct points. We say $C$ is of banana type if it is constructed from a banana curve by attaching at least one rational curve. Similarly, we say a stable map $f:C\to X$ is a banana curve (resp. curve of banana type) if $C$ is a banana curve (resp. curve of banana type). }

{\prop[Nodal Bend-and-Break Lemma]  \label{nodalMBB lemma} 
Let $T \subset \overline{M}^{bir}_{0,0}(X,e)$ be the closure of the $1$-dimensional locus of stable maps of nodal rational type going through $e-1$ general points for $e\geqslant 3$. Suppose $T$ parametrizes a stable map $f: Z \to X$ such that $Z$ contains at least two irreducible components $Z_0$ and $Z_1$. If $f|_{Z_i}$'s are free, then $Z = Z_0 + Z_1$, $f$ is birational onto its image, and either
\begin{enumerate}
    \item $f|_{Z_0}$ is a stable map of nodal rational type, or
    
    \item $f$ factors through a banana curve.
\end{enumerate}
}

{\proof 
Write $f: Z = \sum_{i=0}^r Z_i + \sum_{r+1}^k Z_i + Z_C$, where $Z_i$'s are free for $i = 0,\dots, r$ and non-free for $i = r+1,\dots,k$, and $Z_C$ are the contracted components. Let $e_i$ denote the $H$-degree of $Z_i$, and $e'$ the total $H$-degree of the non-free components.

Since $Z$ goes through $e-1$ general points and each $Z_i$ goes through at most $e_i$ general points, we have the following inequality:
\[
e - 1 \leqslant \sum_{i=0}^r e_i = e-e'.
\]
This implies that $e' \leqslant 1$ and there is at most one non-free line. If there is such a non-free line $L$, then each free rational curve goes through the maximum number of general points. Hence there are finitely many choices for each $Z_i$ and they are disjoint from each other. Since there is a $1$-dimensional family of non-free lines and it is a codimension one condition for each member of this family to intersect the $Z_i$'s, there exist finitely many choices for $L$ upon meeting one of the $Z_i$'s. However, this implies that the image of $Z$ is not connected because there are at least two free components. Hence there is no non-free component and all but one of the free rational components must go through $e_i$ general points. The remaining component, say $Z_0$, must contain $e_0-1$ general points, and hence deforms in a $2$-dimensional family. 

Suppose $Z$ contains at least two other components $Z_1$ and $Z_2$ beside $Z_0$. Then intersecting $Z_1$ and $Z_2$ imposes a two dimensional constraint on the deformation of $Z_0$. Since the points we chose are general, the image of $Z$ has genus zero. This is a contradiction to the assumption that $f$ is a stable map of nodal rational type. Therefore $Z$ has exactly two components. If $f|_{Z_0}$ is a stable map of nodal rational type, then $Z_0$ must goes through at least two general points, and hence it deforms in a $1$-dimensional family by Proposition \ref{free genus one nodal stable map} and Lemma \ref{dimension of nodal rational curve}. Thus it intersects $Z_1$ at a single point. If the images of $Z_0$ and $Z_1$ are smooth rational curves, then by Lemma \ref{tangent intersection has high codimension}, $Z_0$ does not intersect $Z_1$ tangentially. Hence they intersect in at least two distinct points, in which case $f$ factors through a banana curve. Finally, a similar argument in the proof of Proposition \ref{Prop:MBB} shows that $f$ is birational onto its image.
\qed}

\smallskip
To further deform the resulting nodal curve, we apply a similar construction as in \cite[Section 5]{LT19a} and \cite{Tes05}. 

\smallskip
{\construction Let $Z = \sum_{i=0}^r Z_i$ be a connected nodal curve of genus $g$ such that $Z_i$ are irreducible components. Let $f: Z\to X$ be a stable map in $\overline{M}_{g,0}(X)$. 

Let $g_i$ be the genus of $Z_i$ and $e_i = \deg f|_{Z_i}$. We can choose irreducible components $M_i \subset \overline{M}_{g_i,0}(X, e_i)$ containing $f|_{Z_i}$. Denote by $M_i^{(n)}$ the induced component in $\overline{M}_{g_i,n}(X)$ of $M_i$. If $M_i$ is a free (resp. very free) component, let $U_i \subset M_i$ denote the open locus parametrizing the locus of free (resp. very free) curves, and $U_i^{(n)}$ the n-pointed locus in $\overline{M}_{g_i,n}(X)$. }

{\definition A connected nodal curve $f:Z\to X$ with $k$ marked points is a chain of length $r+1$ if the following conditions hold:
\begin{itemize}
    \item The domain $Z = \sum_{i=0}^r$ is a union of $r$  irreducible components $Z_i$.
    \item For $i\neq j$, $Z_i\cap Z_j$ is either empty or a single point.
    \item The map $f$ can be parametrized by some point in the locus of the form 
    \[
    M_0^{(k+1)}\times_X M_1^{(2)} \times_X \cdots \times_X M_{r-1}^{(2)} \times_X M_r^{(1)},
    \]
    where $M_j^{(n_j)} \subset \overline{M}_{g_j,n_j}(X)$. In particular, all $k$ marked points lie on $Z_0$ and $f|_{Z_0} \in M_0^{(k)}$.
\end{itemize}
An end of a chain is an irreducible component that intersects the rest in exactly one point.} 

{\definition Let  
\[
M \subset M_0^{(k+1)}\times_X M_1^{(2)} \times_X \cdots \times_X M_{r-1}^{(2)} \times_X M_r^{(1)}
\]
be an irreducible component such that the projection map from $M$ to each factor in the product is dominant. We call $M$ a principal component of the product. }

\smallskip
Let $f:Z = \sum_{i=0}^r Z_i \to X$ be a chain of length $r+1$. Suppose $f|_{Z_i}$'s are free for all $i$. We may decorate $f$ with marked points all lying on an end $Z_0$ of $Z$. Then $[f]$ can be regarded as a smooth point in a principal component 
\[
U \subset U_0^{(k+1)} \times_X U_1^{(2)} \times_X \dots \times_X U_r^{(1)},
\]
where $k$ is the number of marked points on $Z_0$. Given such a stable map, we are interested in the following deformation:

{\definition[Sliding move] \label{Slide} Let $[g]$ be a stable map. We say $[g]$ is obtained from $[f]$ by sliding $Z_0$ along $Z_i$ if:
\begin{enumerate}[\hspace{2mm}1.]
    \item The domain $Z'$ of $[g]$ can be written as union of connected components $Z' = Z_0' + \sum_{i=1}^r Z_i'$ such that $Z_i'\cong Z_i$ for all $i$.
    \item The maps $g|_{Z_i'} \cong f|_{Z_i}$ for $i = 1,\dots, r$.
    \item The map $[g]$ is contained in $U$.
\end{enumerate}}

\smallskip
Additionally, when dealing with deformations of banana curves, we also need to slide $Z_0$ along $Z_i$ while maintaining $Z_0$ going through some fixed point in $X$. More precisely, let $U_0(p)$ denote the irreducible component in $U_0$ parametrizing deformations of $[f]$ going through a fixed point $f(p)$ for some $p\in Z_0$. Let $U(p)$ denote the fibre product of $U_0(p)^{(k+1)}\times_{U_0^{(k+1)}} U$, i.e, the sublocus of $U$ parametrizing chains going through $f(p)$.

{\definition[Sliding move with point fixed] We say $[g]$ is obtained from $[f]$ by sliding $Z_0$ along $Z_i$ with a point fixed if conditions 1. and 2. in Definition \ref{Slide} hold, and in addition, the map $[g]$ is contained in $U(p)$.}

\smallskip
The following lemma tells us when can we perform these sliding moves:

{\lemma \label{slidingmove} Let $f: Z \to X$ be a chain parametrized by some principal component $M$. Let $Z_0$ denote an end of $Z$ and let $Z_1$ be the component intersecting $Z_0$. If $f|_{Z_0}$ is free, then we can slide $Z_0$ along $Z_1$. If  $f|_{Z_0}$ is very free, then given a smooth point $p\in Z$ that lies on $Z_0$, we can slide $Z_0$ along $Z_1$ with $f(p)$ fixed.}

{\proof Both claims are a consequence of the fact that the locus of geometrically irreducible fibres is a constructible set. In particular, for the second claim, let $U_0(p)$ denote the irreducible locus containing deformations of $[f]$ with $f(p)$ fixed. Since $f|_{Z_0}$ is very free, $U_0(p)$ generically parametrizes maps with irreducible domain and sweeps out a dominant family on $X$, hence we can slide $Z_0$ along $Z_1$ with $f(p)$ fixed.
\qed}

\smallskip
It is also useful to interchange the components in the domain of a chain for smoothing purposes. The following definition makes this precise:

{\definition \cite[c.f. Definition 5.10]{LT19a} Let $f: Z = Z_1 + \dots Z_r \to X$ be a chain. Suppose $f|_{Z_i}$ is contained in $\overline{M}_{g_i, n_i}(X)$ for each $i$. Let $f^\dagger: \{1,\dots, r\} \to \{\overline{M}_{g_1, n_1}(X), \dots, \overline{M}_{g_r, n_r}(X)\}$ be the function which assigns to $i$ the unique irreducible component in $\overline{M}_{g_i, n_i}(X)$ containing $f|_{Z_i}$. We call $f^\dagger$ the combinatorial type of $f$.
}

\smallskip
Specifically, we may exchange the components in a chain of free curves. The proof of the following lemma is similar to that of \cite[Lemma 5.11]{LT19a}, and we omit it.

{\lemma \cite[Lemma 5.11]{LT19a} \label{exchange combinatorial type} Let $M\subset \overline{M}_{1,0}(X)$ be an irreducible component. Suppose $M$ contains a chain $f: Z = Z_1 + \dots + Z_r \to X$ of free curves. Then for any $\tilde{f}^\dagger$ which is a precomposition of $f^\dagger$ with a permutation, $M$ also parametrizes a chain of free curves with combinatorial type $\tilde{f}^\dagger$. }

\smallskip
The next lemma shows that a banana curve is a smooth point of the moduli space if the restricted tangent bundle to each component is positive enough.

{\lemma \label{smoothbanana} Let $f: Z = Z_1 + Z_2 \to X$ be a banana curve. Suppose $f|_{Z_i}$ are free for $i=1,2$ and $f|_{Z_i}$ is very free for at least one $Z_i$. Then $f$ is a smooth point in $\overline{M}_{1,0}(X)$.}

{\proof Without loss of generality, assume $Z_1$ is very free and $Z_2$ is free. Let $p$ and $q$ denote the point of intersections. Tensoring the short exact sequence 
\[
0 \to \CO_{Z_1}(-p-q) \to \CO_Z \to \CO_{Z_2} \to 0
\]
by $f^*T_X$ yields
\[
0 \to f^*T_X|_{Z_1}(-p-q) \to f^*T_X \to f^*T_X|_{Z_2} \to 0.
\]
By assumption, $H^1(Z_1, f^*T_X|_{Z_1}(-p-q)) = 0$ and $H^1(Z_2, f^*T_X|_{Z_2}) = 0$, so $H^1(Z, f^*T_X) = 0$, and $f$ is a smooth point in $\overline{M}_{1,0}(X)$.
\qed}

\smallskip
We are ready to deform our banana curves. In particular, we partially normalize the domain of a banana curve so that we may view the new map as a stable map in $\overline{M}_{0,0}(X)$. We perform the sliding moves in $\overline{M}_{0,0}(X)$ and lift the resulting map back to $\overline{M}_{1,0}(X)$. We need the following lemma that guarantees our deformation stays in the same irreducible component.

{\lemma \label{smoothbtcontraction} Let $f: Z\to X$ be a curve of banana type. Suppose $f$ contracts a connected component $R$ of genus zero. Let $Z_i$ be the connected components of $Z' = \overline{Z\backslash R}$, for $i = 0,\dots, m$. Assume that each $Z_i$ intersects $R$ at exactly one point for $i\geqslant 1$. If the restriction of $f$ to each irreducible component of the $Z_i$'s is free and the restriction of $f$ to each irreducible component of $Z_0$ is very free, then $f$ is a smooth point in $\overline{M}_{1,0}(X)$.}

{\proof Let $P$ denote the subscheme of $Z'$ representing the intersection of $R$ with the other components in $Z$. We have the exact sequence
\[
0 \to f^*T_X|_{Z'}(-P) \to f^*T_X \to f^*T_X|_{R} \to 0.
\]
Since $R$ is a rational curve, we have $h^1(R, f^*T_X|_{R}) = 0$. Since $Z'$ is the disjoint union of the $Z_i$'s, we have
\[
h^1(Z', f^*T_X|_{Z'}(-P)) = \sum_{i=0}^m h^1(Z', f^*T_X|_{Z_i}(-P_i)).
\]
where $P_i = P|_{Z_i}$. By assumption, $P_i$ has length one for $i = 1,\dots, m$, and since $f$ is a genus one banana curve, $P_0$ has length at most $2$. Hence $h^1(Z', f^*T_X|_{Z_i}(-P_i)) = 0$ for all $i$.
\qed}

{\prop \label{bananasmoothie} Let $f: Z = Z_1 + Z_2 \to X$ be a banana curve of degree $e\geqslant 5$. Let $p$ and $q$ be the intersection points of $Z_1$ and $Z_2$. Suppose that $f$ is an immersion and that $p$ and $q$ are mapped to general points in $X$. Then $f$ deforms to a curve of banana type $g$ satisfying the following conditions:
\begin{itemize}
    \item The domain of $g$ contains three irreducible components $Z_1'$, $Z_2'$, and $Z_2$, where $Z_2$ and $Z_2'$ form a banana curve. 

    \item $g$ is an immersion.
    
    \item $[g]$ is a smooth point in the same component of $\overline{M}_{1,0}(X)$ containing $[f]$.
\end{itemize}
}

{\proof Let $f_i$ denote the restriction of $f$ to $Z_i$. Since $e\geqslant 5$, we may assume $f_1$ has degree at least $3$. Apply \cite[Proposition 6.7]{BLRT22} to $f_1$ with the points $f(p)$ and $f(q)$ fixed. We obtain a broken curve $f': Z' = Z'_1 + Z'_2 + Z_2 \to X$, where $Z_1$ deforms to $Z'_1 + Z'_2$, and $f'|_{Z'_1}$ and $f'|_{Z'_2}$ are free curves. Furthermore, we have $p = Z'_1\cap Z_2$ and $q = Z'_2 \cap Z_2$. We also denote $s = Z'_1 \cap Z'_2$. The curve $f'$ lies in the same component as $f$ by \cite[Lemma 5.9]{LT19a}.

Next, we perform the sliding moves by regarding $f'$ as a stable map in $\overline{M}_{0,0}(X)$. Since $f_1$ has degree at least $3$, we may assume $f'|_{Z_2'}$ has degree at least $2$, and hence is very free. We normalize the domain $Z'$ of $f'$ at the intersection point $q$ so that it becomes a genus zero curve $\overline{Z}'$. Let $\overline{Z}_2$ and $\overline{Z}'_2$ be the preimages of $Z_2$ and $Z'_2$ under the normalization. Let $q_1\in \overline{Z}_2$ and $q_2\in \overline{Z}'_2$ be the preimages of $q$. We denote this stable map as $\tilde{f}':\overline{Z}' \to X$. 

By Lemma \ref{slidingmove}, we can slide $\overline{Z}'_2$ along $Z'_1$ with $q$ fixed until $s$ reaches the point $p$, i.e. the images of $Z'_1$, $\overline{Z}'_2$, and $\overline{Z}_2$ in $X$ intersect at the same point $f(p)$. Let $\tilde{f}''$ denote this resulting stable map. Then the domain of $\tilde{f}''$ consists of a rational curve $R$ attached to $Z_1'$, $\overline{Z}_2$, and $\overline{Z}'_2$ at $r_1$, $r_2$, and $r_2'$ respectively. Furthermore, $\tilde{f}''$ contracts $R$. We claim that $\tilde{f}''|_{\overline{Z}_2+\overline{Z}_2'}$ is an immersion. Let $e_2 = \tilde{f}''|_{\overline{Z}_2}$ and $e_2' = \tilde{f}''|_{\overline{Z}_2}$. The claim is clear if $e_2 \neq e_2'$ or $e_2 = e_2' \geqslant 3$. When $e_2 = e_2' = 2$, there are $6$ conics through two general points by \cite[Section 3]{Bea95}. Hence we may slide $\overline{Z}_2'$ so that $\tilde{f}''|_{\overline{Z}_2'}$ has a different image from $\tilde{f}''|_{\overline{Z}_2}$. 

On the other hand, let $\tilde{g}: \overline{Z}' = Z'_1 + \overline{Z}_2 + \overline{Z}'_2 \to X$ be a stable map of genus zero such that $p = \overline{Z}_2\cap \overline{Z}'_2$, $s = Z'_1\cap Z_2'$, and $\tilde{g}|_{\overline{Z}_2+\overline{Z}_2'}$ has the same image as $\tilde{f}''|_{\overline{Z}_2+\overline{Z}_2'}$. Then $\tilde{f}''$ is also obtained as the limit of $\tilde{g}$ by sliding $\overline{Z}'_2$ along $\overline{Z}_2$ until $s$ reaches $p$.

\begin{figure}[h]
    \centering
    \begin{minipage}{0.32\textwidth}
        \centering
        \begin{tikzpicture}[scale=0.8]
            \node[circle, fill=black, inner sep=2pt, label=below left:$Z_1'$] (Z1p) at (0,0) {};
            \node[circle, fill=black, inner sep=2pt, label=below right:$Z_2$] (Z2) at (3,0) {};
            \node[circle, fill=black, inner sep=2pt, label=above:$Z_2'$] (Z2p) at (1.5,2) {};
        
            \draw (Z1p) -- node[below] {$p$} (Z2);
            \draw (Z2) -- node[right] {$q$} (Z2p);
            \draw (Z1p) -- node[left] {$s$} (Z2p);
        \end{tikzpicture}
    \end{minipage}
    \hfill
    \begin{minipage}{0.32\textwidth}
        \centering
        \begin{tikzpicture}[scale=1.25]
            \node[circle, fill=black, inner sep=2pt, label=above left:$Z_1'$] (Z1p) at (0,1) {};
            \node[circle, fill=black, inner sep=2pt, label=below left:$Z_2$] (Z2) at (0,-1) {};
            \node[circle, fill=black, inner sep=2pt, label=right:$Z_2'$] (Z2p) at ({sqrt(3)},0) {};
        
            \node[circle, fill=black, inner sep=2pt, label=above:$R$] (R) at ({sqrt(3)/3},0) {};
        
            \draw (R) -- node[left] {$r_1$} (Z1p);
            \draw (R) -- node[left] {$r_2$} (Z2);
            \draw (R) -- node[above right] {$r_2'$} (Z2p);
        
            \draw (Z2p) -- node[below] {$q$} (Z2);
        \end{tikzpicture}
    \end{minipage}
    \hfill
    \begin{minipage}{0.32\textwidth}
        \centering
        \begin{tikzpicture}[scale=0.8]
        \node[circle, fill=black, inner sep=2pt, label=above:$Z_2'$] (Z2p) at (0,1.5) {};
        \node[circle, fill=black, inner sep=2pt, label=below left:$Z_1'$] (Z1p) at (-1.5,-1) {};
        \node[circle, fill=black, inner sep=2pt, label=below right:$Z_2$] (Z2) at (1.5,-1) {};
        
        \draw (Z2p) to[bend left=25] node[right] {$p$} (Z2);
        \draw (Z2p) to[bend right=25] node[left] {$q$} (Z2);
        
        \draw (Z1p) -- node[left] {$s$} (Z2p);
        \end{tikzpicture}
    \end{minipage}

    \caption{Dual graphs of $f'$, $f''$, and $g$ respectively.}
\end{figure}

Finally, let $f''$ and $g$ denote the corresponding stable maps of $\tilde{f}''$ and $\tilde{g}$ in $\overline{M}_{1,0}(X)$ by gluing $\overline{Z}_2$ and $\overline{Z}'_2$ at the $q_i$'s. Then $g$ is a curve of banana type. By Lemma \ref{smoothbanana} and Lemma \ref{smoothbtcontraction}, the stable maps $f''$ and $g$ are smooth points in $\overline{M}_{1,0}(X)$. Hence we conclude that the deformation from $f$ to $g$ is contained in the same irreducible component in $\overline{M}_{1,0}(X)$.
\qed}

{\theorem[Nodal Bend-and-Break]\label{nodalMainMBB} Let $f: C \to X$ be a general nodal rational curve of degree $e\geqslant 5$. Then $f$ deforms to a stable map $g: Z = Z_1 + Z_2 \to X$ that is a smooth point of $\widetilde{M}_{1,0}^{bir}(X)$ satisfying the following conditions:
\begin{itemize}
    \item $Z_1$ is a smooth genus one curve and $g|_{Z_1}$ is free.
    \item $Z_2$ is a rational curve and $g|_{Z_2}$ is free.
    \item $g|_{Z_i}$'s are birational onto their images.
\end{itemize}
}

{\proof By \cite[Lemma 4.1]{LT23} and Proposition \ref{nodalMBB lemma}, the stable map $f$ deforms in $\overline{M}_{0,0}(X)$ to either $f_1: Z = Z_1 + Z_2 \to X$ where $f|_{Z_1}$ is a stable map of nodal rational type, or a banana curve $f_2: Z = Z_1 + Z_2 \to X$. In particular, we can lift both $f_1$ and $f_2$ to $\overline{M}_{1,0}(X)$ by lifting $f_1$ to a stable map $f'_1: C + Z_2 \to X$, where $C$ is an irreducible nodal rational curve of genus one, and lifting $f_2$ to a stable map $f_2': B \to X$, where $B$ is a banana curve, respectively.

The argument to both cases are similar, and we only prove the second case. Since $[f_2']$ is a smooth point in $\overline{M}_{1,0}(X)$ by Lemma \ref{smoothbanana}, we may further deform $f_2'$ so that the nodes of the banana curve are mapped to general points of $X$. Then by Proposition \ref{bananasmoothie}, $f_2'$ deforms to a curve of banana type $f_2'': Z = Z_1 + Z_2 + Z_3 \to X$ such that $[f_2''|_{Z_2 + Z_3}]$ is a smooth point of $\overline{M}_{1,0}(X)$. Since $f_2''$ is an immersion, we can smooth $f_2''|_{Z_2 + Z_3}$ to a stable map sending a smooth genus one curve birationally onto its image. This produces the desired $g$. 
\qed}

\smallskip
We are ready to prove the Movable Bend-and-Break Theorem \ref{Thm:MBB}.

{\proof[Proof of Theorem \ref{Thm:MBB}] 
Let $M\subset \widetilde{M}^{bir}_{1,0}(X,e)$ be a free component with $e\geqslant 5$. Let $g: E\to X$ be a general member of $M$. Let $T$ be the locus containing $[g]$ in the statement of Proposition \ref{Prop:MBB}. By Corollary \ref{Cor: mor does not dominate kontsevich}, the component $M$ generically parametrizes free curves, hence $h^1(E,g^*T_X) = 0$ and the map $H^0(E, N_g) \to H^1(E,T_E)$ is a surjection. This implies that the image of the universal family $\CU_T$ of $T$ under the projection $\pi: \CU_M \to \overline{M}_{1,1}$ is equal to $\overline{M}_{1,1}$. Then $T$ either parametrizes a stable map with a reducible domain which is precisely the desired broken curve by Proposition \ref{Prop:MBB}, or parametrizes a stable map $f$ whose image under $\pi$ is the nodal rational curve $[C, p]\in \overline{M}_{1,1}$. In the latter case, a dimension count shows that the locus of maps with reducible domain above $[C, p]$ has codimension at least $2$. Since $T$ is a $1$-dimensional family constructed from general points and a general basepoint free curve, the domain of $f$ is $C$. Moreover, $[f]$ is a smooth point of $M$ by Proposition \ref{free genus one nodal stable map}. Hence we may apply nodal Bend-and-Break \ref{nodalMainMBB} to $[f]$ to obtain the desired broken curve.
\qed}

\subsection{Proof of main theorems}\hfill

\smallskip
\noindent We are ready to prove our main theorems. The following lemma will be useful in our induction step:

{\lemma \cite[Lemma 2.1]{LT21} \label{usefullemmaflat} Let $M_1, M_2$, and $X$ be irreducible schemes of finite type  with flat and dominant morphisms $f_1: M_1\to X$ and $f_2: M_2\to X$. Let $M\subset M_1\times_X M_2$ be an irreducible component. Suppose a general fibre of $f_2$ is purely $k$-dimensional. Then the natural map $M\to X$ is dominant and $M$ has dimension $k+\dim M_1$. Furthermore, if a general fibre of $f_2$ is irreducible, then $M = M_1\times_X M_2$. }

{\proof[Proof of Theorem \ref{IrredM10}] We prove by induction on the degree. The base case is Corollary \ref{base1}. Suppose now that $\widetilde{M}^{bir}_{1,0}(X,e')$ is irreducible for every degree $e' < e$.

Let $M \subset \widetilde{M}^{bir}_{1,0}(X,e)$ be an irreducible component. Let $[f]\in M$ be general. Movable Bend-and-Break \ref{Thm:MBB} implies that $f$ deforms to a stable map $g: Z_1 + Z_2 \to X$, where $Z_1$ is a smooth genus one curve, $Z_2$ is a rational curve, and the restrictions $g|_{Z_i}$ are free curves. If the degree of $g|_{Z_1}$ is greater than $5$, than by induction on lower degree, we can deform $Z_1$ into a reducible chain $Z_1' + Z_2'$ while fixing the intersection with $Z_2$. Without loss of generality, assume $Z_1'$ is the smooth genus one curve, and by Lemma \ref{exchange combinatorial type}, we may further assume that $Z_1'$ is an end of the chain. Since $\widetilde{M}^{bir}_{0,0}(X,e')$ is irreducible by Theorem \ref{Starr}, we can smooth $Z_2 + Z_2'$ to an irreducible rational curve, which we again denote by $Z_2'$. Hence we can deform $g: Z_1 + Z_2 \to X$ to $g': Z_1' + Z_2' \to X$ with the property that $Z_1$ and $Z_1'$ are smooth genus one curves and $\deg g|_{Z_1} > \deg g'|_{Z_1'}$.

We may repeatedly performing the above procedure until we arrive at a stable map $h: \widetilde{Z_1} + \widetilde{Z_2} \to X$ such that $[h|_{\widetilde{Z_1}}] \in \widetilde{M}^{bir}_{1,0}(X,4)$ and $[h|_{\widetilde{Z_2}}] \in \widetilde{M}^{bir}_{0,0}(X,e - 4)$. Note that if the broken genus one component $Z_1'$ has degree less than $4$, then we may break the rational component $Z_2'$ into a chain of lines and smooth $Z_1'$ and the line attached to it into a stable map $[h|_{\widetilde{Z_1}}] \in \widetilde{M}^{bir}_{1,0}(X,4)$. In particular, $[h|_{\widetilde{Z_1}}]$, and $[h|_{\widetilde{Z_2}}]$ are smooth points in their corresponding components. 

By \cite[Theorem 7.9]{LT19a}, a general fibre of the evaluation map for $\widetilde{M}^{bir, \circ}_{0,1}(X,e - 4)$ is irreducible. Hence by Lemma \ref{usefullemmaflat}, $[h]$ is a smooth point in the irreducible locus $N = \widetilde{M}^{bir, \circ}_{1,1}(X,4) \times_X \widetilde{M}^{bir, \circ}_{0,1}(X,e - 4)$, where $\widetilde{M}^{bir, \circ}_{g,1}(X)$ denotes the locus of free curves with a marked point. Since any irreducible component $M \subset \widetilde{M}^{bir}_{1,0}(X,e)$ is a free component by Theorem \ref{non-free-main}, and a general point in $M$ deforms to a smooth point in the irreducible locus $N$, we conclude that $\widetilde{M}^{bir}_{1,0}(X,e)$ must be irreducible.
\qed}

{\proof[Proof of Theorem \ref{IrredMor}] By Theorem \ref{IrredM10}, there is a unique irreducible component $R_e \subset \widetilde{M}_{1,0}(X,e)$ generically parametrizing free curves of degree $e\geqslant 5$. 

Let $\CR_e \subset \overline{\CM}_{1,0}(X,e)$ be the stack whose coarse moduli space is $R_e$. Let $\CR'_e \subset \widetilde{\CM}_{1,1}(X,e)$ denote its universal family. The composition of the forgetful morphism $\overline{\CM}_{1,1}(X,e) \to \overline{\CM}_{1,1}$ and the coarse moduli space map $\overline{\CM}_{1,1} \to \overline{M}_{1,1}$ gives a morphism
\[
\pi: \CR'_e \to \overline{M}_{1,1}.
\]
At the level of geometric points, $\pi$ sends a stable map $f:Z\to X$ to a point $(E,p)$ in $\overline{M}_{1,1}$, where $E$ represents the isomorphism class of the stabilization of $Z$. We can identify an open substack of the fibre above a geometric point $[E,p] \in \overline{M}_{1,1}$ with the locus $M \subset \Mor(E,X)$ parametrizing degree $e$ maps that are birational onto their images. Upon replacing $\CR'_e$ by a resolution of singularity $\widetilde{\CR}'_e$, we obtain a morphism $\widetilde{\pi}: \widetilde{\CR}'_e \to \overline{M}_{1,1}$. The irreducibility of a general fibre of $\pi$ then follows from the connectedness of a general fibre of $\widetilde{\pi}$.

Suppose to the contrary that a general fibre of $\widetilde{\pi}$ is disconnected. Then the Stein factorization for Deligne-Mumford stacks \cite[Theorem 4.6.14]{Alp25} produces a morphism $\phi: \widetilde{\CR}'_e \to Y$ and a morphism $\psi: Y \to \overline{M}_{1,1}$, where $\phi$ has connected fibres and $\psi$ is a finite morphism of degree greater than one. Let $B \subset \overline{M}_{1,1}$ be the branch locus of $\psi$. 

Let $[E,p]$ be a general point in $\overline{M}_{1,1}$ and let $f: E \to X$ be a stable map representing a general point in the fibre of $\widetilde{\pi}$ above $[E,p]$. By Proposition \ref{bal}, the normal bundle $N_f$ is balanced. Hence we can construct a $1$-dimensional locus $T\subset \widetilde{M}^{bir}_{1,0}(X,e)$ parametrizing stable maps passing through $e-1$ general points and a basepoint free curve as in the statement of Proposition \ref{Prop:MBB}. We may regard $T$ as a closed substack of $\widetilde{\CR}'_e$. Since $\psi$ is finite of degree greater than one, there exists a everywhere non-reduced component $N$ in the fibre $F_b = \widetilde{\CR}'_e \times_{\overline{M}_{1,1}} \{b\}$ above some branch point $b\in B$ such that $T\cap N$ is nonempty. Let $g:Z\to X$ be the stable map representing a point in $T\cap N$. Since the support of $B$ consists of at least two distinct points, we may assume that the stabilization of $Z$ is a smooth genus one curve. By Proposition \ref{Prop:MBB}, either $Z$ is irreducible, or $Z = Z_0 + Z_1$ such that $g|_{Z_0}$ is a free smooth genus one curve and $g|_{Z_1}$ is a free rational curve. We divide our discussion based on the possibilities of $Z$.

If $Z$ is irreducible, then $N$ is a non-free component parametrizing maps of degree at least $6$. This contradicts Theorem \ref{no non-free-comp-morphism}. If $Z$ is reducible, then $N$ generically parametrizes stable maps with reducible domain. Let $g_0 = g|_{Z_0}$ and $g_1 = g|_{Z_1}$. Since the $g_i$'s are free, the exact sequence
\[
0 \to g_1^*T_X(-p) \to g^*T_X \to g_0^*T_X \to 0
\]
shows that $H^1(Z, g^*T_X) = H^1(Z_0, g_0^*T_X) = 0$, contradicting that $N$ is everywhere non-reduced.

\qed}

{\remark The number of irreducible components is upper semi-continuous in a family. The above argument gives the irreducibility of a general fibre, and it would be interesting to show the statement for all fibres.}

\bibliography{ref}
\bibliographystyle{alphaurl}

\end{document}